\documentclass[12pt,a4paper,reqno]{amsart} 
\usepackage{amssymb}
\usepackage{latexsym}
\usepackage{amsmath}
\usepackage{mathrsfs}
\usepackage{cite}

\usepackage{calc}                   

\newcommand{\modBr}[3]{M^{#1}(\omega _{#3},\mathscr{B}^{#2})}
\newcommand{\modBrk}[3]{M^{#1}(\omega _{#3},{#2},\mathscr{B})}
\newcommand{\back}[1]{\setminus #1}
\newcommand{\Char}{\operatorname{Char}}
\newcommand{\scal}[2]{\langle #1,#2\rangle}
\newcommand{\rr}[1]{\mathbf R^{#1}}

\newcommand{\nm}[2]{\Vert #1\Vert _{#2}}

\newcommand{\op}{\operatorname{Op}}

\newcommand{\sets}[2]{\{ \, #1\, ;\, #2\, \} }

\newcommand{\fy}{\varphi}
\newcommand{\cdo}{\, \cdot \, }
\newcommand{\supp}{\operatorname{supp}}

\newcommand{\eabs}[1]{\langle #1\rangle}     

\newcommand{\vrum}{\vspace{0.1cm}}

\newcommand{\ttbigcup}{{\textstyle{\, \bigcup \, }}}
\newcommand{\FB}{\mathscr F\! \mathscr B}
\newcommand{\WF}{\operatorname{WF}}
\newcommand{\dir}{{\operatorname{dir}}}

\numberwithin{equation}{section}          

\newtheorem{thm}{Theorem}
\numberwithin{thm}{section}
\newtheorem*{tom}{\rubrik}
\newcommand{\rubrik}{}
\newtheorem{prop}[thm]{Proposition}

\newtheorem{lemma}[thm]{Lemma}

\theoremstyle{definition}

\newtheorem{defn}[thm]{Definition}

\theoremstyle{remark}

\newtheorem{rem}[thm]{Remark}              

\author{Sandro Coriasco}

\address{Department of mathematics, Turin's University, Italy}

\email{sandro.coriasco@unito.it}

\author{Karoline Johansson}

\address{Department of Mathematics and Systems Engineering,
V{\"a}xj{\"o} University, Sweden}

\email{karoline.johansson@vxu.se}

\author{Joachim Toft}

\address{Department of Mathematics and Systems Engineering,
V{\"a}xj{\"o} University, Sweden}

\email{joachim.toft@vxu.se}

\title{Wave-front sets of Banach function types}



\begin{document}

\begin{abstract}
Let $\omega ,\omega _0$ be appropriate weight functions and $\mathscr B$ be
an invarian BF-space.
We introduce the wave-front set, $\WF_{\FB(\omega )}(f)$ of the distribution $f$
with respect to weighted Fourier Banach space $\FB(\omega )$. We prove
that usual mapping properties for pseudo-differential operators
$\op _t (a)$ with symbols $a$ in $S^{(\omega _0)}_{\rho ,0}$ hold for
such wave-front sets. In particular we prove
$\displaystyle
\WF _{\FB(\omega /\omega _0)}(\op _t (a) f)\subseteq \WF _{\FB(\omega
)}(f)$ and $\WF _{\FB(\omega
)}(f)
\subseteq \WF _{\FB(\omega /\omega
_0)}(\op _t (a) f)\ttbigcup \Char (a)$.
Here $\Char (a)$ is the set of characteristic points of $a$.
\end{abstract}

\maketitle

\section{Introduction}\label{sec0}

\par

In this paper we introduce
wave-front sets with respect to Fourier images of translation
invariant BF-spaces. The family of such
wave-front sets contains the wave-front sets of Sobolev type,
introduced by H{\"o}rmander in \cite {Hrm-nonlin}, the
classical wave-front sets (cf. Sections 8.1
and 8.2 in \cite {Ho1}), and wave-front sets of Fourier Lebesgue
types, introduced in \cite{PTT1}. Roughly speaking, for any
given distribution $f$ and for appropriate Banach (or Frech\'et) space
$\mathcal B$ of tempered distributions, the wave-front set
$\WF  _{\mathcal B}(f)$ of $f$ consists of all pairs
$(x_0,\xi _0)$ in $\rr d\times (\rr d\back 0)$ such that no
localizations of the distribution at $x_0$ belongs to $\mathcal B$ in
the direction $\xi _0$.

\par

We also establish mapping properties for a quite general class of
pseudo-differential operators on such wave-front sets, 
and show that the micro-local analysis in \cite{PTT1} in background
of Fourier Lebesgue spaces can be further generalized. It follows that our approach
gives rise to flexible micro-local analysis tools which fit well to
the most common approach developed in e.g. \cite {Ho1,
Hrm-nonlin}. In particular, we prove that
usual mapping properties, which are valid for classical wave-front
sets (cf. Chapters VIII and XVIII in \cite{Ho1}), also hold for
wave-front sets of Fourier Banach types.
For example, we show
\begin{equation}\label{eq:inclusions}
\begin{aligned}
\WF _{\FB(\omega /\omega _0)}(\op _t (a) &f)\subseteq \WF _{\FB(\omega
)}(f)\\
&
\subseteq \WF _{\FB(\omega /\omega
_0)}(\op _t (a) f)\ttbigcup \Char (a).
\end{aligned}
\end{equation}
That is, any operator
$\op (a)$ shrinks the wave-front sets and opposite embeddings can be
obtained by including $\Char (a)$, the set of characteristic points of
the operator symbol $a$.

\par

The symbol classes for the pseudo-differential operators are denoted by
$S_{\rho ,\delta}^{(\omega _0)}(\rr {2d})$, the set of all
smooth functions $a$ on $\rr {2d}$ such that $a/\omega _0\in S^0_{\rho
,\delta}(\rr {2d})$. Here $\rho ,\delta \in \mathbf R$ and $\omega_0$ is
an appropriate smooth function on $\rr
{2d}$. We note that $S_{\rho ,\delta}^{(\omega _0)}(\rr {2d})$ agrees
with the H{\"o}rmander class $S^r_{\rho ,\delta}(\rr {2d})$ when
$\omega _0(x, \xi )=\eabs \xi ^r$, where $r\in \mathbf R$ and $\eabs \xi
=(1+|\xi |^2)^{1/2}$.

\par

The set of characteristic points $\Char (a)$ of $a\in
S_{\rho ,\delta}^{(\omega )}$ is the same as in \cite{PTT1}, and
depends on the choices of
$\rho$, $\delta$ and $\omega$ (see Definition \ref{defchar} and
Proposition \ref{psiecharequiv}).
We recall that this set is smaller than the set of characteristic
points given by \cite{Ho1}. It is empty when $a$ satisfies a local
ellipticity condition
with respect to $\omega$, which is fulfilled for any hypoelliptic
partial differential operator with constant coefficients
(cf. \cite{PTT1}). As a consequence of \eqref{eq:inclusions}, it follows that such hypoelliptic
operators preserve the wave-front sets, as expected (cf. Example
3.9 in \cite{PTT1}).

\medspace

Information on regularity in background of wave-front sets of
Fourier Banach types might be more detailed compared to classical
wave-front sets, because of our choices of different weight functions
$\omega$ and Banach spaces when defining our Fourier Banach space
$\FB(\omega )(\rr d)$. For example, the space $\FB(\omega)=\mathscr
FL^1_{(\omega )}(\rr d)$, with $\omega (x,\xi )=\eabs \xi ^N$ for some
integer $N\ge 0$, is locally close to $C^N(\rr d)$ (cf. the Introduction of
\cite{PTT1}). Consequently, the wave-front set with respect to $\mathscr FL^1_{(\omega )}$ can be
used to investigate a sort of regularity which is close to smoothness
of order $N$.

\medspace

Furthermore, we are able to apply our results on pseudo-differential
operators in context of modulation space theory, when discussing
mapping properties of pseudo-differential operators with respect to
wave-front sets.  
The modulation spaces were introduced by Feichtinger in \cite{F1}, and
the theory was developed in
\cite{Feichtinger3, Feichtinger4, Feichtinger5, Grochenig0a}. The
modulation space $M(\omega, \mathscr B )$, where $\omega$ is
a weight function (or time-frequency shift) on phase space $\rr {2d}$,
appears as the set of temperated (ultra-)distributions
whose short-time Fourier transform belong to the weighted Banach space
$\mathscr B(\omega )$.
These types of modulation spaces contains the (classical) modulation
spaces $M^{p,q}_{(\omega)}(\rr {2d})$ as well as the space
$W^{p,q}_{(\omega)}(\rr {2d})$ related to the Wiener amalgam spaces,
by choosing $\mathscr B =L^{p,q}_1(\rr {2d})$ and $\mathscr B
=L^{p,q}_2(\rr {2d})$ respectively (see Remark \ref{Modamalgam}).
In the last part of the paper we define wave-front sets with respect to
weighted modulation spaces, and prove that they coincide with the
wave-front sets of Fourier Banach types.

\par

Parallel to this development, modulation spaces have been
incorporated into the calculus of pseudo-differential operators, in
the sense of the study of continuity of (classical)
pseudo-differential operators acting on modulation spaces (cf.
\cite{Tachizawa1,Czaja,Pilipovic2,Pilipovic3,Teofanov1,Teofanov2}),
and the study of operators of non-classical type, where modulation
spaces are used as symbol classes. We refer to \cite{Gro-book, Grochenig2, Grochenig0,
Grochenig1b, Grochenig1c, Herau1, HTW, Pilipovic2,
Sjostrand1, Sjostrand2, Toft2, Toft35, To8, Toft4} for more facts about pseudo-differential operators in background of modulation space theory.

\medspace

The paper is organized as follows. In Section
\ref{sec1} we recall the definition and basic properties for
pseudo-differential operators, translation invariant Banach function
spaces (BF-spaces) and (weighted) Fourier Banach
spaces. Here we also define sets of characteristic points for a broad
class of pseudo-differential operators.
In Section \ref{sec2} we prove some properties for the sets of
characteristic points, which shows that our definition coincide with
the sets of characteristic points defined in \cite{PTT1}. These sets
might be smaller than characteristic sets in
\cite{Ho1} (cf. \cite[Example 3.11]{PTT1}).

\par 

In Section \ref{sec3} we define wave-front sets with respect to
(weighted) Fourier Banach spaces, and prove some of their
main properties. Thereafter, in Section \ref{sec4} we
show how these wave-front sets are propagated under the action
of pseudo-differential operators. In particular, we prove \eqref{eq:inclusions},
when $\omega _0$ and $\omega$ are appropriate weights and
$a$ belongs to $S^{(\omega _0)}_{\rho ,0}$ with
$\rho >0$.

\par

In Section \ref{sec5} we consider wave-front sets obtained from
sequences of Fourier Banach spaces. These types of wave-front sets
contain the classical ones (with respect to smoothness), and the mapping
properties for pseudo-differential operators also hold in this context
(cf. Section 18.1 in \cite{Ho1}).

\par

Finally, Section \ref{sec6} is devoted to study the definition and basic properties of
wave-front sets with respect to modulation spaces. We prove that they can be identified with certain
wave-front sets of Fourier Banach types.

\par

\section{Preliminaries}\label{sec1}

\par

In this section we recall some notation and basic results. The proofs
are in general omitted. In what
follows we let $\Gamma$ denote an open cone in $\rr d\back 0$. If
$\xi \in \rr d\back 0$ is fixed, then an open cone which contains
$\xi $ is sometimes denoted by $\Gamma_\xi$.

\par

Assume that $\omega, v\in L^\infty _{loc}(\rr d)$ are positive
functions. Then $\omega$ is called $v$-moderate if
\begin{equation}\label{moderate}
\omega (x+y) \leq C\omega (x)v(y)
\end{equation}
for some constant $C$ which is independent of $x,y\in \rr d$. If $v$
in \eqref{moderate} can be chosen as a polynomial, then $\omega$ is
called polynomially moderate. We let $\mathscr P(\rr d)$ be the set
of all polynomially moderated functions on $\rr d$. We say that $v$ is
\emph{submultiplicative} when \eqref{moderate} holds with $\omega =v$.
Throughout we assume that the submultiplicative weights are even.
If $\omega (x,\xi )\in \mathscr P(\rr {2d})$ is constant with respect
to the $x$-variable ($\xi$-variable), then we sometimes write $\omega
(\xi )$ ($\omega (x)$) instead of $\omega (x,\xi )$. In this case we
consider $\omega$ as an element in $\mathscr P(\rr {2d})$ or in
$\mathscr P(\rr d)$ depending on the situation.

\par

We also need to consider classes of weight functions, related to
$\mathscr P$. More precisely, we let $\mathscr P_0(\rr d)$ be the set
of all $\omega \in \mathscr P(\rr d)\bigcap C^\infty (\rr d)$ such
that $\partial ^\alpha \omega /\omega \in L^\infty$ for all
multi-indices $\alpha$.
For each $\omega \in \mathscr P(\rr d)$, there is an equivalent
weight $\omega _0\in \mathscr P_0(\rr d)$, that is,
$C^{-1}\omega _0\le \omega \le C\omega _0$ holds for some constant
$C$ (cf. \cite[Lemma 1.2]{To8}).

\par

Assume that $\rho ,\delta \in \mathbf R$. Then we let $\mathscr
P_{\rho ,\delta}(\rr {2d})$ be the set of all $\omega (x,\xi )$ in
$\mathscr P(\rr {2d})\cap C^\infty (\rr {2d})$ such that
$$
\eabs \xi ^{\rho |\beta |-\delta |\alpha |}\frac {(\partial ^\alpha
_x\partial ^\beta _\xi \omega )(x,\xi )}{\omega (x,\xi )}\in L^\infty
(\rr {2d}),
$$
for every multi-indices $\alpha$ and $\beta$. Note that in contrast to
$\mathscr P_0$, we do not have an equivalence between $\mathscr
P_{\rho ,\delta}$ and $\mathscr P$ when $\rho >0$. On the other hand,
if $s\in \mathbf R$ and $\rho \in [0,1]$, then $\mathscr P_{\rho
,\delta} (\rr {2d})$ contains $\omega (x,\xi )=\eabs \xi ^s$, which
are one of the most important classes in the applications.

\par

For any weight $\omega$ in $\mathscr P(\rr d)$ or in $\mathscr P_{\rho
,\delta}(\rr d)$, we let $L^p_{(\omega )}(\rr d)$ be the set of all
$f\in L^1_{loc}(\rr d)$ such that $f\cdot \omega \in L^p(\rr d)$.

\medspace

The Fourier transform $\mathscr F$ is the linear and continuous
mapping on $\mathscr S'(\rr d)$ which takes the form
$$
(\mathscr Ff)(\xi )= \widehat f(\xi ) \equiv (2\pi )^{-d/2}\int _{\rr
{d}} f(x)e^{-i\scal  x\xi }\, dx
$$
when $f\in L^1(\rr d)$. We recall that $\mathscr F$ is a homeomorphism
on $\mathscr S'(\rr d)$ which restricts to a homeomorphism on $\mathscr
S(\rr d)$ and to a unitary operator on $L^2(\rr d)$.

\medspace

Next we recall the definition of Banach function spaces.

\par

\begin{defn}\label{BFspaces}
Assume that $\mathscr B$ is a Banach space of complex-valued
measurable functions on $\rr d$ and that $v \in \mathscr P(\rr {d})$
is submultiplicative.
Then $\mathscr B$ is called a \emph{(translation) invariant
BF-space on $\rr d$} (with respect to $v$), if there is a constant $C$
such that the following conditions are fulfilled:
\begin{enumerate}
\item $\mathscr S(\rr d)\subseteq \mathscr
B\subseteq \mathscr S'(\rr d)$ (continuous embeddings);

\vrum

\item if $x\in \rr d$ and $f\in \mathscr B$, then $f(\cdot -x)\in
\mathscr B$, and
\begin{equation}\label{translmultprop1}
\nm {f(\cdot -x)}{\mathscr B}\le Cv(x)\nm {f}{\mathscr B}\text ;
\end{equation}

\vrum

\item if $f,g\in L^1_{loc}(\rr d)$ satisfy $g\in \mathscr B$  and $|f|
\le |g|$ almost everywhere, then $f\in \mathscr B$ and
$$
\nm f{\mathscr B}\le C\nm g{\mathscr B}\text .
$$
\end{enumerate}
\end{defn}

\par

Assume that $\mathscr B$ is a translation invariant BF-space. If $f\in
\mathscr B$ and $h\in L^\infty$, then it follows from (3) in
Definition \ref{BFspaces} that $f\cdot h\in \mathscr B$ and
\begin{equation}\label{multprop}
\nm {f\cdot h}{\mathscr B}\le C\nm f{\mathscr B}\nm h{L^\infty}.
\end{equation}

\par

\begin{rem}\label{newbfspaces}
Assume that $\omega _0,v,v_0\in \mathscr P(\rr d)$ are such $v$ and
$v_0$ are submultiplicative,
$\omega _0$ is $v_0$-moderate, and assume that $\mathscr B$ is a
translation-invariant BF-space on $\rr d$ with respect to $v$. Also
let $\mathscr B_0$ be the Banach space which consists of all $f\in
L^1_{loc}(\rr d)$ such that $\nm f{\mathscr B_0}\equiv \nm {f\, \omega _0
}{\mathscr B}$ is finite. Then $\mathscr B_0$ is a translation
invariant BF-space with respect to $v_0v$.
\end{rem}

\par

\begin{rem}\label{BFemb}
Let $\mathscr B$ be an invariant BF-space. Then it is easy to find
Sobolev type spaces which are continuously embedded in $\mathscr
B$. In fact, for each $p\in [1,\infty ]$ and integer $N\ge 0$, let
$Q^p_N(\rr d)$ be the set of all $f\in L^p(\rr d)$ such that $\nm
f{Q^p_N}<\infty$, where
$$
\nm f{Q^p_N}\equiv \sum _{|\alpha +\beta |\le N}\nm {x^\alpha
D^\beta f}{L^p}.
$$
Then for each $p$ fixed, the topology for $\mathscr S(\rr d)$ can be
defined by the semi-norms $f\mapsto \nm f{Q^p_N}$, for $N=0,1,\dots
$.

\par

A combination of this fact and (1) and (3) in Definition
\ref{BFspaces} now shows that for each $p\in [1,\infty ]$ and each
translation invariant BF-space $\mathscr B$, there is an integer $N\ge
0$ such that $Q^p_N(\rr d)\subseteq \mathscr B$. Moreover, let
$L^\infty _N(\rr d)$ be the set of all $f\in L^\infty _{loc}(\rr d)$
such that $f\, \eabs \cdo ^N\in L^\infty$. Then, since any
element in $L^\infty _N$ can be majorized with an element in
$Q^\infty _N$, it follows from (3) in Definition \ref{BFspaces} that
$L^\infty _N\subseteq \mathscr B$, provided $N$ is chosen large
enough. This proves the assertion.
\end{rem}

\par

For future references we note that if $\mathscr B$ is a translation
invariant BF-space with respect to the submultiplicative weight $v$ on
$\rr d$, then the convolution map $*$ on $\mathscr S(\rr d)$ extends
uniquely to a continuous mapping from $\mathscr B\times L^1_{(v)}(\rr
d)$, and for some constant $C$ it holds
\begin{equation}\label{propupps}
\nm {\fy *f}{\mathscr B}\le C\nm \fy{L^1_{(v)}}\nm f{\mathscr B},\quad
\fy \in L^1_{(v)}(\rr d),\ f\in \mathscr B.
\end{equation}
In fact, if $f\in \mathscr B$ and $g$ is a step function, then $f*g$ is
well-defined and belongs to $\mathscr B$ in view of the definitions, and
Minkowski's inequality gives
\begin{multline*}
\nm {f*g}{\mathscr B} =\Big \Vert \int f(\cdo -y)g(y)\, dy \Big \Vert
_{\mathscr B}
\\[1ex]
\le \int \nm { f(\cdo-y)}{\mathscr B}|g(y)|\, dy \le C\int \nm {
f}{\mathscr B}|g(y)v(y)|\, dy = C\nm { f}{\mathscr B}\nm
g{L^1_{(v)}}.
\end{multline*}

\par

Now assume that $g\in C_0^\infty$. Then $f*g$ is well-defined as an
element in $\mathscr S' \cap C^\infty$, and by approximating $g$ with
step functions and using \eqref{propupps} it follows that $f*g\in
\mathscr B$ and that \eqref{propupps} holds also in this case.
The assertion now follows from this fact and a simple argument of
approximations, using the fact that $C^\infty _0$ is dense in
$L^1_{(v)}$.

\par

For each translation invariant BF-space $\mathscr B$ on $\rr d$, and
each pair of vector spaces $(V_1,V_2)$ such that $V_1\oplus V_2=\rr
d$, we define the projection spaces $\mathscr B_1$ and
$\mathscr B_2$ of $\mathscr B$ by the formulae
\begin{alignat}{2}
&\mathscr{B}_1 &  &\equiv \sets {f\in \mathscr
S'(V_1)}{f\otimes \fy \in \mathscr B \text{ for every } \fy \in 
\mathscr S (V_2)}\label{B1def}
\intertext{and}
&\mathscr{B}_2 &  &\equiv \sets {f\in \mathscr
S'(V_2)}{\fy \otimes f \in \mathscr B \text{ for every } \fy \in 
\mathscr S (V_1)}.\label{B2def}
\end{alignat}

\par

\begin{prop}\label{propbnoll}
Let $\mathscr{B}$ be a translation
invariant BF-space on $\rr d$, and let $\mathscr B_1$ and $\mathscr
B_2$ be the same as in \eqref{B1def} and \eqref{B2def}. Then 
\begin{alignat}{2}
&\mathscr{B}_1 &  &= \sets {f\in \mathscr
S'(V_1)}{f\otimes \fy \in \mathscr B \text{ for some } \fy \in 
\mathscr S (V_2)\back 0}\tag*{(\ref{B1def})$'$}
\intertext{and}
&\mathscr{B}_2 &  &= \sets {f\in \mathscr
S'(V_2)}{\fy \otimes f \in \mathscr B \text{ for some } \fy \in 
\mathscr S (V_1)\back 0}.\tag*{(\ref{B2def})$'$}
\end{alignat}

\par

In particular, if $\fy _j\in \mathscr S (V_j)\back 0$ for $j=1,2$ are
fixed and $f_1\in \mathscr{S}'(V_1)$ and $f\in \mathscr{S}'(V_2)$,
then $\mathscr B_1$ and $\mathscr B_2$ are translation
invariant BF-spaces under the norms
$$
\nm f{\mathscr B_1}\equiv \nm {f\otimes \fy _1}{\mathscr B}\quad
\text{and}\quad \nm f{\mathscr B_2}\equiv \nm {\fy _2\otimes
f}{\mathscr B}
$$
respectively.
\end{prop}

\par

\begin{proof}
We only prove \eqref{B2def}$'$. The other equality follows by similar
arguments and is left for the reader. We may assume that $V_j=\rr
{d_j}$ with $d_1+d_2=d$.

\par

Let $\mathscr B_0$ be the right-hand side of \eqref{B2def}$'$. Then it
is obvious that $\mathscr B_2\subseteq \mathscr B_0$. We have to prove
the opposite inclusion.

\par

Therefore, assume that $f\in \mathscr B_0$, and choose $\fy_0 \in
\mathscr S(\rr{d_1})\back 0$ such that $\fy _0\otimes f \in
\mathscr{B}$. Also let $\fy \in \mathscr S(\rr{d_1})$ be
arbitrary. We shall prove that $\fy \otimes f \in
\mathscr{B}$.

\par

Let $Q \subseteq \rr {d_1}$ be an open ball and $c>0$
be chosen such that $|\fy _0(x)| >c$ when $x\in Q$. Also
let the lattice $\Lambda \subseteq \rr {d_1}$ and $\fy _1\in
C^{\infty}_0(Q)$ be such that $0\leq \fy _1\leq 1$ and
$$
\sum_{\{x_j\}\in \Lambda} \fy _1(\cdo -x_j)=1.
$$
Then $\fy _1\leq C|\fy _0|$, for some constant $C>0$, which gives
$$
\nm {\fy _1\otimes f}{\mathscr B}\le C\nm {\fy
_0\otimes f}{\mathscr B}<\infty .
$$
This in turn gives
\begin{multline}\label{chi0est1}
\nm {\fy \otimes f}{\mathscr B} \le \sum \nm {(\fy _1(\cdo -x_j)\fy
)\otimes f}{\mathscr B}
\\[1ex]
\le \sum  v(x_j,0)\nm {(\fy _1\fy (\cdo +x_j))\otimes f}{\mathscr B} 
\\[1ex]
\le C
\Big (\sum  v(x_j,0)\nm {\fy (\cdo +x_j)}{L^\infty (Q)}\Big ) \nm
{\fy _1\otimes f}{\mathscr B}.
\end{multline}
Since $v\in \mathscr P$ and $\fy \in \mathscr S$, it follows that the
sum in the right-hand side of \eqref{chi0est1} is finite. Hence $f\in
\mathscr B_2$, and the proof is complete.
\end{proof}

\par

\begin{rem}
We note that the last sum in \eqref{chi0est1} is the norm
$$
\nm \fy {W_{(v)}} \equiv
\sum  v(x_j,0)\nm {\fy (\cdo +x_j)}{L^\infty (Q)}
$$ 
for the weighted Wiener space
$$
W_{(v)}(\rr d) =\sets {f\in L^\infty _{loc}(\rr d)}{\nm f{W_{(v)}}<\infty }
$$
(cf. \cite{Gro-book}). The results in Proposition \ref{propbnoll} can therefore be improved in such way that we may replace $\mathscr S$ by $W_{(v)}$ in \eqref{B1def}, \eqref{B2def}, \eqref{B1def}$'$ and \eqref{B2def}$'$.
\end{rem}

\par

Assume that $\mathscr B$ is a translation invariant BF-space on $\rr
d$, and that $\omega \in \mathscr P(\rr d)$. Then we let $\FB
{(\omega )}$ be the set of all $f\in \mathscr S'(\rr d)$ such that
$\xi \mapsto \widehat f(\xi )\omega (x,\xi )$ belongs to $\mathscr
B$. It follows that $\FB {(\omega )}$ is a Banach space under the
norm
\begin{equation}\label{FLnorm}
\nm f{\FB {(\omega )}}\equiv \nm {\widehat
f\, \omega}{\mathscr B}.
\end{equation}

\par

\begin{rem} \label{whyomega}
In many situations it is convenient to permit an $x$
dependency for the weight $\omega$ in the definition of Fourier Banach spaces.
More precisely, for each $\omega \in \mathscr P(\rr
{2d})$ and each translation invariant BF-space $\mathscr B$ on $\rr d$, we let $\mathscr {FB}{(\omega )}$ be the set of all $f\in
\mathscr S'(\rr d)$ such that
$$
\nm f{\mathscr {FB}{(\omega )}}
= \nm f{\mathscr {FB}{(\omega),x}}
\equiv \nm {\widehat f\, \omega (x,\cdo )}{\mathscr B}
$$
is finite. Since $\omega$ is
$v$-moderate for some $v\in \mathscr P(\rr {2d})$ it follows that
different choices of $x$ give rise to equivalent norms.
Therefore the condition
$\nm f{\FB(\omega )}<\infty$ is
independent of $x$, and  it follows that $\FB(\omega )(\rr
d)$ is  independent of $x$ although $\nm \cdo {\FB(\omega )}$ might depend on $x$.
\end{rem}

\par

Recall that a topological vector space $V\subseteq \mathscr
D'(X)$ is called \emph{local} if $V\subseteq V_{loc}$. Here
$X\subseteq \rr d$ is open, and $V_{loc}$ consists of all $f\in
\mathscr D'(X)$ such that $\fy f \in V$ for every $\fy \in C_0^\infty
(X)$. For future references we note that if $\mathscr B$ is a
translation invariant BF-space on $\rr d$ and $\omega \in \mathscr
P(\rr {2d})$, then it follows from \eqref{propupps} that $\FB(\omega
)$ is a local space, i.{\,}e.
\begin{equation}\label{Blocal}
\FB (\omega )\subseteq \FB (\omega )_{loc}\equiv (\FB (\omega ))_{loc}.
\end{equation}

\medspace

We need to recall some facts from Chapter XVIII in \cite {Ho1}
concerning pseudo-differential operators. Let $a\in
\mathscr S(\rr {2d})$, and $t\in \mathbf R$ be fixed. Then
the pseudo-differential operator $\op _t(a)$ is the linear and
continuous operator on $\mathscr S(\rr d)$, defined by the formula
\begin{equation}\label{e0.5}
(\op _t(a)f)(x)
=
(2\pi ) ^{-d}\iint a((1-t)x+ty,\xi )f(y)e^{i\scal {x-y}\xi }\,
dyd\xi .
\end{equation}
For general $a\in \mathscr S'(\rr {2d})$, the pseudo-differential
operator $\op _t(a)$ is defined as the continuous operator from
$\mathscr S(\rr d)$ to $\mathscr S'(\rr d)$ with distribution
kernel
\begin{equation}\label{weylkernel}
K_{t,a}(x,y)=(2\pi )^{-d/2}(\mathscr F_2^{-1}a)((1-t)x+ty,x-y).
\end{equation}
Here $\mathscr F_2F$ is the partial Fourier transform of
$F(x,y)\in \mathscr S'(\rr{2d})$ with respect to the $y$-variable.
This definition makes sense, since the mappings $\mathscr F_2$ and
$$
F(x,y)\mapsto F((1-t)x+ty,x-y)
$$
are homeomorphisms on $\mathscr
S'(\rr {2d})$. We also note that the latter definition of $\op _t(a)$
agrees with the operator in \eqref{e0.5} when $a\in \mathscr S(\rr
{2d})$. If  $t=0$, then $\op _t(a)$ agrees with the Kohn-Nirenberg
representation $\op (a)=a(x,D)$.

\par

If $a\in \mathscr S'(\rr {2d})$ and $s,t\in
\mathbf R$, then there is a unique $b\in \mathscr S'(\rr {2d})$ such
that $\op _s(a)=\op _t(b)$. By straight-forward applications of
Fourier's inversion  formula, it follows that
\begin{equation}\label{pseudorelation}
\op _s(a)=\op _t(b) \quad \Longleftrightarrow \quad b(x,\xi
)=e^{i(t-s)\scal
{D_x}{D_\xi}}a(x,\xi ).
\end{equation}
(Cf. Section 18.5 in \cite{Ho1}.)

\par

Next we discuss symbol classes which we use. Let $r,
\rho ,\delta \in \mathbf R$ be fixed. Then we recall from
\cite{Ho1} that $S^r_{\rho ,\delta}(\rr {2d})$ is the set of all $a\in
C^\infty (\rr {2d})$ such that for each pairs of multi-indices
$\alpha$ and $\beta$, there is a constant $C_{\alpha ,\beta}$ such
that
$$
|\partial _x^\alpha \partial _\xi ^\beta a(x,\xi )|\le
C_{\alpha ,\beta }\eabs \xi ^{r-\rho |\beta |+\delta |\alpha |}.
$$
Usually we assume that $0\le \delta \le \rho \le 1$, $0<\rho$ and
$\delta <1$.

\par

More generally, assume that $\omega \in \mathscr P_{\rho ,\delta } (\rr
{2d})$. Then we recall from the introduction that
$S_{\rho ,\delta}^{(\omega )}(\rr {2d})$ consists of all $a\in
C^\infty (\rr {2d})$ such that
\begin{equation}\label{Somegadef}
|\partial _x^\alpha \partial _\xi ^\beta a(x,\xi )|\le
C_{\alpha ,\beta }\omega (x,\xi )\eabs \xi ^{-\rho |\beta
|+\delta |\alpha |}.
\end{equation}
We note that $S_{\rho ,\delta}^{(\omega )}(\rr {2d})=S(\omega ,g_{\rho
,\delta})$, when $g=g_{\rho 
,\delta}$ is the Riemannian metric on $\rr
{2d}$, defined by the formula
$$
\big (g_{\rho ,\delta }\big )_{(y,\eta )}(x,\xi ) = \eabs \eta
^{2\delta}|x|^2 + \eabs \eta ^{-2\rho}|\xi |^2
$$
(cf. Section 18.4--18.6 in \cite{Ho1}). Furthermore, $S^{(\omega
)}_{\rho ,\delta} =S^r_{\rho ,\delta}$ when $\omega (x,\xi )=\eabs \xi
^r$, as remarked in the introduction.

\par

The following result shows that pseudo-differential operators with
symbols in $S^{(\omega )}_{\rho ,\delta}$ behave well. We refer to
\cite {Ho1} or \cite{PTT1} for the proof.

\par

\begin{prop}\label{p5.4}
Let $\rho ,\delta \in [0,1]$ be such that $0\le \delta \le
\rho \le 1$ and $\delta <1$, and let $\omega \in \mathscr P_{\rho ,\delta}
(\rr {2d})$. If $a\in S^{(\omega )}_{\rho ,\delta} (\rr
{2d})$, then $\op _t(a)$ is continuous on $\mathscr S(\rr d)$ and
extends uniquely to a continuous operator on $\mathscr S'(\rr d)$.
\end{prop}

\par

We also need to define the set of characteristic points of a symbol
$a\in S^{(\omega )}_{\rho ,\delta} (\rr {2d})$, when $\omega \in
\mathscr P_{\rho ,\delta}(\rr {2d})$ and $0\le \delta < \rho \le
1$. In Section \ref{sec2} we show that this definition is equivalent
to Definition 1.3 in \cite{PTT1}. We remark that our sets of
characteristic points are smaller than the corresponding sets in \cite
{Ho1}. (Cf. \cite[Definition 18.1.5]{Ho1} and Remark
\ref{compchar} in Section \ref{sec2}).

\par

\begin{defn}\label{defchar}
Assume that  $0 \leq \delta < \rho \leq 1$, $\omega _0 \in
\mathscr{P}_{\rho,\delta}(\rr{2d})$ and $a\in
S^{(\omega_0)}_{\rho,\delta}(\rr{2d})$. Then $a$ is called
\emph{$\psi$-invertible} with respect to $\omega _0$ at the
point $(x_0,\xi_0)\in \rr d\times (\rr d\back 0)$, if
there exist a neighbourhood $X$ of $x_0$, an open conical
neighbourhood $\Gamma$ of $\xi_0$ and positive constants $R$ and $C$ such that
\begin{equation}\label{nydef}
|a(x,\xi)|\geq C\omega _0(x,\xi),
\end{equation}
for $x\in X$, $\xi\in \Gamma$ and $|\xi|\geq R$.

\par

The point $(x_0,\xi_0)$ is called \emph{characteristic} for $a$ with
respect to $\omega _0$ if $a$ is \emph{not} $\psi$-invertible 
with respect to $\omega _0$ at $(x_0,\xi_0)$. The set of characteristic points (the
characteristic set), for $a$ with respect to $\omega _0$ is denoted
$\Char(a)=\Char_{(\omega _0)}(a)$.
\end{defn}

\par

We note that $(x_0,\xi _0)\notin \Char_{(\omega _0)}(a)$ means that $a$
is elliptic near $x_0$ in the direction $\xi _0$. Since the case $\omega _0=1$ in Definition \ref{defchar} is especially important we also make the following definition. We say that $c\in
S^{0}_{\rho,\delta}(\rr{2d})$ is
\emph{$\psi$-invertible} at $(x_0,\xi_0) \in \rr d\times (\rr d\back 0)$, if $(x_0,\xi _0)\notin \Char _{(\omega _0)}(c)$ with $\omega _0=1$. That is,  there exist a neighbourhood $X$ of $x_0$, an open conical
neighbourhood $\Gamma$ of $\xi_0$ and $R >0$ such that \eqref{nydef} holds for $a=c$ and $\omega _0=1$,
for some constant $C>0$ which is independent of $x\in X$
and $\xi\in \Gamma$ such that $|\xi|\geq R$.

\par

It will also be convenient to have the following definition of different types of cutoff functions.

\par

\begin{defn}\label{cuttdef}
Let $X\subseteq \rr d$ be open, $\Gamma \subseteq \rr d\back 0$ be an open cone, $x_0\in X$ and let $\xi _0\in \Gamma $. 

\begin{enumerate}
\item A smooth function $\fy$ on $\rr d$ is called a cutoff function with respect to $x_0$ and $X$, if $0\le \fy \le 1$, $\fy \in C_0^\infty (X)$ and $\fy =1$ in an open neighbourhood of $x_0$. The set of cutoff functions with respect to $x_0$ and $X$ is denoted by $\mathscr C_{x_0}(X)$;

\vrum

\item  A smooth function $\psi$ on $\rr d$ is called a directional cutoff function with respect to $\xi_0$ and $\Gamma$, if there is a constant $R>0$ and open conical neighbourhood $\Gamma _1$ of $\xi _0$ such that the following is true:
\begin{itemize}
\item $0\le \psi \le 1$ and $\supp \psi \subseteq \Gamma$;

\vrum

\item  $\psi (t\xi )=\psi (\xi )$ when $t\ge 1$ and $|\xi |\ge R$;

\vrum

\item $\psi (\xi )=1$ when $\xi \in \Gamma _1$ and $|\xi |\ge R$.
\end{itemize}

The set of directional cutoff functions with respect to $\xi _0$ and $\Gamma$ is denoted by $\mathscr C^\dir  _{\xi _0}(\Gamma )$.
\end{enumerate}

\end{defn}

\par

\begin{rem}\label{psiinvremark}
We note that if $\fy \in \mathscr C_{x_0}(X)$ and $\psi \in \mathscr
C^\dir  _{\xi _0}(\Gamma )$ for some $(x_0,\xi _0)\in \rr d\times (\rr
d\back 0)$, then $c\equiv \fy \otimes \psi$ belongs to $S^0_{1,0}(\rr
{2d})$ and is $\psi$-invertible at $(x_0,\xi _0)$.
\end{rem}

\par

\section{Pseudo-differential calculus with symbols
in $S^{(\omega )}_{\rho ,\delta}$}\label{sec2}

\par

In this section we make a review of basic results for
pseudo-differential operators with symbols in classes of the form
$S^{(\omega )}_{\rho ,\delta}(\rr {2d})$, when $0\le \delta <\rho \le
1$ and $\omega \in \mathscr P_{\rho ,\delta}(\rr {2d})$. For the
standard properties in the pseudo-differential calculus we only state
the results and refer to \cite{Ho1} for the proofs. Though there are
similar stated and
proved properties concerning sets of characteristic points, we include
proofs of these properties in order to being more self-contained.

\par

We start with the following result concerning compositions and
invariance properties for pseudo-differential operators. Here we let
$$
\sigma _s(x,\xi )=\eabs \xi ^s,
$$
where $\eabs \xi =(1+|\xi
|^2)^{1/2}$ as usual. We also recall that $S^{-\infty }_{\rho ,\delta
}=S^{-\infty }_{1,0}$ consists of all $a\in C^\infty (\rr {2d})$ such
that for each $N\in \mathbf R$ and multi-index $\alpha$, there is a
constant $C_{N,\alpha}$ such that
$$
|\partial ^\alpha a(x,\xi )|\le C_{N,\alpha}\eabs \xi ^{-N}.
$$

\par

\begin{prop}\label{pseudocomp}
Let  $0\le \delta <\rho \le 1$, $\mu =\rho -\delta >0$ and
$\omega ,\omega _1,\omega _2\in \mathscr P_{\rho ,\delta }(\rr
{2d})$. Also let $\{ m_j\} _{j=0}^{\infty}$ be a sequence of real
numbers such that $m_j\to -\infty$ as $j\to \infty$. Then the
following is true:
\begin{enumerate}
\item if $a_1\in S^{(\omega _1)}_{\rho ,\delta }(\rr {2d})$ and
$a_2\in S^{(\omega _2)}_{\rho ,\delta }(\rr {2d})$, then $\op
(a_1)\circ \op (a_2)=\op (c)$, for some $c\in S^{(\omega
_1\omega _2)}_{\rho ,\delta }(\rr {2d})$. Furthermore,
\begin{equation}\label{symbcomp}
c(x,\xi )-\sum _{|\alpha |<N}\frac {i^{|\alpha |}(D^{\alpha}_\xi
a_1)(x,\xi )(D^{\alpha}_x a_2)(x,\xi )}{\alpha
!}\in S^{(\omega
_1\omega _2\sigma _{-N\mu })}_{\rho ,\delta }(\rr {2d})
\end{equation}
for every $N\ge0$;

\vrum

\item if $M=\sup _{k\ge 0}(m_k)$,
$M_j=\sup _{k\ge j}(m_k)$ and $a_j\in S^{(\omega \sigma _{m_j})}_{\rho
,\delta }(\rr {2d})$, then it exists $a\in S^{(\omega \sigma
_{M})}_{\rho ,\delta }(\rr {2d})$ such that
\begin{equation}\label{asymptexp}
a(x,\xi )-\sum _{|\alpha |<N}a_j(x,\xi )\in S^{(\omega 
\sigma _{M_N})}_{\rho ,\delta }(\rr {2d});
\end{equation}
for every $N\ge 0$;

\vrum

\item if $a,b\in \mathscr S'(\rr {2d})$ and $s,t\in \mathbf R$ are
such that $\op _s(a)=\op _t(b)$, then $a\in S^{(\omega )}_{\rho
,\delta }(\rr {2d})$, if and only if $b\in S^{(\omega )}_{\rho ,\delta
}(\rr {2d})$, and
\begin{equation}\label{pseudorel2}
b(x,\xi )-\sum _{k<N}\frac {(i(t-s)\scal {D_x}{D_\xi})^{k}a(x,\xi )}{k
!}\in S^{(\omega \sigma _{-N\mu })}_{\rho ,\delta }(\rr
{2d})
\end{equation}
for every $N\ge 0$.
\end{enumerate}
\end{prop}

\par

As usual we write
\begin{equation}\tag*{(\ref{asymptexp})$'$}
a\sim \sum a_j
\end{equation}
when \eqref{asymptexp} is fulfilled for every $N\ge 0$. In particular
it follows from \eqref{symbcomp} and \eqref{pseudorel2} that
\begin{equation}\tag*{(\ref{symbcomp})$'$}
c\sim \sum \frac {i^{|\alpha |}(D^{\alpha}_\xi
a_1)(D^{\alpha}_x a_2)}{\alpha !}
\end{equation}
when $\op (a_1)\circ \op (a_2)=\op (c)$, and
\begin{equation}\tag*{(\ref{pseudorel2})$'$}
b\sim \sum \frac {(i(t-s)\scal {D_x}{D_\xi})^{k}a}{k !}
\end{equation}
when $\op _s(a)=\op _t(b)$.

\par

In the following proposition we show that the set of characteristic
points for a pseudo-differential operator is independent of the choice
of pseudo-differential calculus.

\par

\begin{prop}\label{anvandningkar}
Assume that $s,t\in \mathbf R$, $0\leq \delta < \rho\leq 1$,
$\omega_0\in\mathscr{P}_{\rho,\delta}$ and that $a,b\in
S^{(\omega_0)}_{\rho,\delta}(\rr{2d})$ satisfy $\op_s(a)=\op_t(b)$. Then
\begin{equation}\label{char}
\Char_{(\omega_0)}(a)=\Char_{(\omega_0)}(b).
\end{equation}
\end{prop}

\par

\begin{proof}
Let $\mu$ and $\sigma _s$ be the same as in the proof of Proposition
\ref{pseudocomp}. By Proposition \ref{pseudocomp} (3) we have
$$
b=a+h,
$$
for some $h\in S^{(\omega_0\sigma_{-\mu})}_{\rho,\delta}$.

\par

Assume that $(x_0,\xi_0)\notin \Char_{(\omega_0)}(a)$. By the
definitions, there is a negihbourhood $X$ of $x_0$, an open conical
negihbourhood $\Gamma$ of $\xi_0$, $C>0$ and $R >0$ such that
$$
|a(x,\xi )| \ge C\omega _0(x,\xi )\quad \text{and}\quad |h(x,\xi )|\le
C\omega _0(x,\xi )/2,
$$
as $x\in X$, $\xi \in \Gamma$ and $|\xi |\ge R$. This gives
$$
|b(x,\xi )| \ge C\omega _0(x,\xi )/2,\quad \text{when}\quad x\in X,\
\xi \in \Gamma ,\ |\xi |\ge R,
$$
and it follows that $(x_0,\xi_0)\notin \Char_{(\omega_0)}(b)$. Hence
$\Char_{(\omega_0)}(b)\subseteq \Char_{(\omega_0)}(a)$. By symmetry, 
the opposite inclusion also holds. Hence
$\Char_{(\omega_0)}(a) = \Char_{(\omega_0)}(b)$, and the proof is
complete.
\end{proof}

\par

The following proposition shows different aspects of set of
characteristic points, and is important when investigating wave-front
properties for pseudo-differential operators. In particular it shows
that $\op (a)$ satisfy certain invertibility properties outside the
set of characteristic points for $a$. More precisely,  
outside $\Char _{(\omega _0)}(a)$, we prove that
\begin{equation}\label{invexp}
\op (b)\op (a) =\op (c)+\op (h),
\end{equation}
for some convenient $b$, $c$ and $h$ which take the role of inverse,
identity symbol and smoothing remainder respectively.

\par

\begin{prop}\label{psiecharequiv}
Let  $0\le \delta <\rho \le 1$, $\omega _0\in \mathscr
P_{\rho,\delta}(\rr {2d})$, $a\in S^{(\omega _0)}_{\rho,\delta}(\rr
{2d})$, $(x_0,\xi _0)\in \rr d \times (\rr d\back 0)$, and let
$\mu =\rho -\delta$. Then the following conditions are equivalent:
\begin{enumerate}
\item $(x_0,\xi _0)\notin \Char _{(\omega _0)} (a)$;

\vrum

\item there is an element $c\in S^0_{\rho ,\delta}$ which is $\psi$-invertible at $(x_0,\xi _0)$, and an element $b\in S^{(1/\omega _0)}_{\rho,\delta}$ such that $ab=c$;

\vrum
 
\item there is an element $c\in S^0_{\rho ,\delta}$ which is $\psi$-invertible at $(x_0,\xi _0)$, and elements 
$h\in S^{-\mu}_{\rho,\delta}$ and $b\in S^{(1/\omega _0)}_{\rho,\delta}$ such that \eqref{invexp} holds;

\vrum

\item for each neighbourhood $X$ of $x_0$ and conical neighbourhood $\Gamma $ of $\xi _0$, there is an element $c=\fy \otimes \psi$ where $\fy \in \mathscr C_{x_0}(X)$ and $\psi _{\xi _0}^\dir (\Gamma )$, and elements 
$h\in \mathscr S$ and $b\in S^{(1/\omega _0)}_{\rho,\delta}$ such that \eqref{invexp} holds. Furthermore, 
the supports of $b$ and $h$ are contained in $X\times \rr d$.
\end{enumerate}
\end{prop}

\par

For the proof we note that $\mu $ in Proposition
\ref{psiecharequiv} is positive, which in turn implies that $\cap _{j\ge 0}S^{(\omega _0\sigma _{-j\mu })}
(\rr {2d})$ agrees with $S^{-\infty}(\rr {2d})$.

\par

\begin{proof}
The equivalence between (1) and (2) follows by letting
$b(x,\xi )=\fy (x)\psi (\xi )/a(x,\xi )$ for some appropriate  $\fy \in \mathscr C_{x_0}(\rr d)$ and $\psi \in \mathscr C_{\xi _0}^\dir (\rr d \back 0)$.

\par

(4) $\Rightarrow$ (3) is obvious in view of Remark \ref{psiinvremark}. Assume that (3) holds. We shall
prove that (1) holds, and since $|b|\leq C/\omega _0$, it suffices to
prove that
\begin{align}
|a(x,\xi )b(x,\xi )| \geq 1/2\label{abinvertible}
\intertext{when}
(x,\xi) \in X\times \Gamma,\ |\xi |\ge R\label{abinvcond}
\end{align}
holds for some conical neighbourhood $\Gamma$ of $\xi _0$, some open
neighbourhood $X$ of $x_0$ and some $R>0$.

\par

By Proposition \ref{pseudocomp} (1) it follows that $ab=c+h$ for
some $h\in S ^{-\mu}_{\rho,\delta}$. By choosing $R$ large enough
and $\Gamma$ sufficiently small conical neighbourhood of $\xi _0$, it
follows that $c(x,\xi )=1$ and $|h(x,\xi )|\leq 1/2$ when
\eqref{abinvcond} holds. This gives \eqref{abinvertible}, and (1)
follows.

\par

It remains to prove that (1) implies (4). Therefore assume that (1) holds, and choose an open neighbourhood $X$ of $x_0$, an open conical neighbourhood $\Gamma $ of $\xi _0$ and $R>0$ such that \eqref{nydef} holds when $(x,\xi )\in X\times \Gamma$ and $|\xi |>R$. Also let $\fy _j\in \mathscr C_{x_0}(X)$ and $\psi _j\in \mathscr C_{\xi _0}^\dir (\Gamma )$ for $j=1,2,3$ be such that $\fy _j=1$ on $\supp \fy _{j+1}$, $\psi _j=1$ on $\supp \psi _{j+1}$ when $j=1,2$, and $\psi _j(\xi )=0$ when $|\xi |\le R$. We also set $c_j=\fy _j\otimes \psi _j$ when $j\le 2$ and $c_j=c_2$ when $j \ge 3$.

\par

If $b_1(x,\xi )=\fy _1(x)\psi _1(\xi )/a(x,\xi
)\in S ^{(1/\omega _0)}_{\rho,\delta}$, then the symbol of $\op
(b_1)\op (a)$ is equal to  $c_1\mod (S ^{-\mu }_{\rho,\delta})$. Hence
\begin{equation}\label{opcompmod}
\op (b_j)\op (a)=\op (c_j)+\op (h_j)
\end{equation}
holds for $j=1$ and some $h_1\in S ^{-\mu }_{\rho,\delta}$.

\par

For $j\geq 2$ we now define $\widetilde b_j\in S ^{(1/\omega _0
)}_{\rho,\delta}$ by the Neumann serie
$$
\op (\widetilde b_j) = \sum _{k=0}^{j-1}(-1)^k\op (\widetilde r_k),
$$
where $\op (\widetilde r_k)=\op (h_1)^k\op (b_1)\in \op (S ^{(\sigma
_{-k\mu}/\omega _0)}_{\rho,\delta})$. Then \eqref{opcompmod} gives
\begin{multline*}
\op (\widetilde b_j)\op (a) = \sum _{k=0}^{j-1}(-1)^k\op (h_1)^k\op
(b_1)\op (a)
\\[1ex]
=\sum _{k=0}^{j-1}(-1)^k\op (h_1)^k(\op (c_1)+\op (h_1)).
\end{multline*}
That is
\begin{equation}\label{expan1}
\op (\widetilde b_j)\op (a) = \op (c_1)+\op (\widetilde h_{1,j})+\op
(\widetilde h_{2,j}),
\end{equation}
where
\begin{align}
\op (\widetilde h_{1,j}) &= (-1)^{j-1}\op (h_1)^{j}\in \op (S
^{-j\mu}_{\rho,\delta})\label{hjtilde}
\intertext{and}
\op (\widetilde h_{2,j}) &= -\sum _{k=1}^{j-1}(-1)^k\op (h_1)^k\op
(1-c_1)\in \op (S^{-\mu}_{\rho,\delta}).\notag
\end{align}

\par

By Proposition \ref{pseudocomp} (1) and asymptotic expansions it
follows that
\begin{multline}\label{h2tilde}
\op (\widetilde h_{2,j}) = -\sum _{k=1}^{j-1}(-1)^k\op
(1-c_1)\op (h_1)^k
\\[1ex]
+ \op (\widetilde h_{3,j}) + \op (\widetilde
h_{4,j}),
\end{multline}
for some $\widetilde h_{3,j}\in S ^{-\mu}_{\rho,\delta}$ which
is equal to zero in $\supp c_1$ and $\widetilde h_{4,j}\in
S^{-j\mu}_{\rho,\delta}$. Now let $b_j$ and $r_k$ be defined by
the formulae
\begin{align*}
\op (b_j) &= \op
(c_2)\op (\widetilde b_j)\in \op (S ^{(1/\omega _0)}_{\rho,\delta})
\\[1ex]
\op (r_k) &= \op (c_2)\op (\widetilde r_k)\in \op (S ^{(\sigma
_{-k\mu }/\omega _0)}_{\rho,\delta}).
\end{align*}
Then
$$
\op (b_j) = \sum _{k=0}^{j-1}(-1)^k\op (r_k)
$$
and \eqref{expan1}--\eqref{h2tilde} give
\begin{multline*}
\op (b_j)\op (a) = \op (c_2)\op (c_1)+\op
(c_2)\op (\widetilde h_{1,j})
\\[1ex]
-\sum _{k=1}^{j-1}(-1)^k\op (c_2)\op
(1-c_1)\op (h_1)^k + \op (c_2)\op (\widetilde h_{3,j}) +
\op (c_2)\op (\widetilde h_{4,j}).
\end{multline*}
Since $c_1=1$ and $\widetilde h_{3,j}=0$ on $\supp c_2$, it
follows that
\begin{align*}
\op (c_2)\op (c_1) &= \op (c_2)\mod \op (S^{-\infty}),
\\[1ex]
\op (c_2)\op (\widetilde h_{1,j}) &\in \op (S^{-j\mu}_{\rho,\delta}),
\\[1ex]
\sum _{k=1}^{j-1}(-1)^k\op (c_2)\op (1-c_1)\op (h_1)^k &\in  \op
(S^{-\infty}),
\\[1ex]
\op (c_2)\op (\widetilde h_{3,j}) &\in \op (S^{-\infty})
\intertext{and}
\op (c_2)\op (\widetilde h_{4,j}) &\in \op (S^{-j\mu}_{\rho,\delta}).
\end{align*}
Hence, \eqref{opcompmod} follows for some $h_j\in
S^{-j\mu}_{\rho,\delta}$.

\par

By choosing $b_0\in S ^{(1/\omega )}_{\rho,\delta}$ such that
$$
b_0\sim \sum r_k,
$$
it follows that $\op (b_0)\op (a) =\op (c_2)+\op (h_0)$, with
$$
h_0\in S ^{-\infty}.
$$
The assertion (4) now follows by letting
\begin{align*}
b(x,\xi )=\fy _3(x)b_0(x,\xi ),\quad  c(x,\xi ) &= \fy _3(x)c_2(x,\xi ),
\\[1ex]
\text{and}\quad h(x,\xi ) &= \fy _3(x) h_0(x,\xi ),
\end{align*}
and using the fact that if $\fy _3 \in C_0^\infty (\rr d)$ and $h_0\in
S^{-\infty}(\rr {2d})$, then $\fy _3 (x)h_0(x,\xi )\in \mathscr S(\rr {2d})$.
The proof is complete.
\end{proof}

\par

\begin{rem}\label{compchar}
By Proposition \ref{psiecharequiv} it follows that Definition 1.3 in
\cite {PTT1} is equivalent to Definition \ref{defchar}. We also remark
that if $a$ is an appropriate symbol, and $\Char '(a)$ the set of
characteristic points for $a$ in the sense of \cite[Definition
18.1.5]{Ho1}, then $\Char _{(\omega _0)}(a)\subseteq \Char
'(a)$. Furthermore, strict embedding might occur, especially for
symbols to hypoelliptic partial operators with constant coefficients, which are not elliptic (cf.  Example 3.11 in \cite{PTT1}).
\end{rem}

\par

\section{Wave front sets with respect to Fourier
Banach spaces}\label{sec3}

\par

In this section we define wave-front sets with respect to Fourier
Banach spaces, and show some basic properties. 

\par

Let $\omega \in \mathscr P(\rr {2d})$, $\Gamma \subseteq \rr
d\back 0$ be an open cone and let $\mathscr B$ be a translation
invariant BF-space on $\rr d$. For any $f\in
\mathscr E'(\rr d)$, let
\begin{equation}\label{skoff}
|f|_{\FB (\omega,\Gamma)}=|f|_{{\FB}(\omega,\Gamma)_x}\equiv \nm
{\widehat{f}\omega(x,\cdo )\chi_{\Gamma} }{\mathscr{B}}.
\end{equation}
We note that $\widehat f\omega(x,\cdo )\chi_{\Gamma}\in \mathscr
B_{loc}$ for every $f\in \mathscr E'$. If $\widehat
f\omega(x,\cdo )\chi_{\Gamma}\notin \mathscr B$, then we set $|f|_{\FB
(\omega,\Gamma)}=+\infty$. Hence  $|\cdo |_{\FB (\omega,\Gamma)}$
defines a semi-norm on $\mathscr E'$ which might attain the value
$+\infty$. Since $\omega $ is $v$-moderate for some $v \in \mathscr
P(\rr {2d})$, it follows that different $x \in \rr d$ gives rise to
equivalent semi-norms. Furthermore, if
$\Gamma =\rr d\back 0$ and $f\in \mathscr {FB}{(\omega )}\cap \mathscr E'$, then
$|f|_{\FB (\omega,\Gamma)}$ agrees with $\nm f{\FB {(\omega)}}$. For simplicity we write $|f|_{\FB (\Gamma)}$ instead of $|f|_{\FB
(\omega,\Gamma)}$ when $\omega =1$.

\par

For the sake of notational convenience we set
\begin{equation} \label{notconv}
|\cdo |_{\mathcal B(\Gamma )}=|\cdo |_{{\FB}(\omega,\Gamma)_x}, \quad
\mbox{when}
\quad
\mathcal B=\FB (\omega).
\end{equation}
We let $ \Theta _{\mathcal B}(f)=\Theta _{{\FB}(\omega)} (f)$ be the set of all $ \xi \in \rr d\back 0 $ such that
$|f|_{\mathcal B(\Gamma )} < \infty$, for some $
\Gamma = \Gamma _{\xi}$.  We also let $\Sigma _{\mathcal B} (f)$
be the complement of $ \Theta_{\mathcal
B} (f)$ in $\rr d\back 0 $. Then
$\Theta _{ \mathcal B} (f)$ and $\Sigma _{\mathcal
B} (f)$ are open respectively
closed subsets in $\rr d\back 0$, which are independent of
the choice of $ x \in \rr d$ in \eqref{skoff}.

\par

\begin{defn}\label{wave-frontdef}
Let  $\mathscr B$ be a translation invariant BF-space on $\rr d$, $\omega \in \mathscr P(\rr {2d})$, $\mathcal
B$ be as in \eqref{notconv}, and let
$X$ be an open subset of $\rr d$.
The wave-front set of
$f\in \mathscr D'(X)$,
$
\WF _{\mathcal B}(f)  \equiv  \WF _{{\FB}(\omega )}(f)
$
with respect to $\mathcal B$ consists of all pairs $(x_0,\xi_0)$ in
$X\times (\rr d \back 0)$ such that
$
\xi _0 \in  \Sigma _{\mathcal B} (\fy f)
$
holds for each $\fy \in C_0^\infty (X)$ such that $\fy (x_0)\neq
0$.
\end{defn}

\par

We note that $\WF  _{\mathcal B}(f)$ in Definition \ref{wave-frontdef}
is a closed set in $X\times
(\rr d\back 0)$, since it is obvious that its complement is
open. We also note that if $ x_0\in \rr d$ is fixed and $\omega _0(\xi
)=\omega (x_0,\xi )$, then $\WF _{{\FB}(\omega )} (f)=\WF _{{\FB}(\omega _0)}(f)$, since $\Sigma _{\mathcal B}$ is independent of $x_0$.

\par

The following theorem shows that wave-front sets with respect to
$\FB(\omega )$ satisfy appropriate micro-local
properties. It also shows that such wave-front sets decreases when the local Fourier BF-spaces increases, or when the weight $\omega$ decreases.

\par

\begin{thm}\label{theta-sigma-propertiesAA}
Let $X\subseteq \rr d$ be open, $\mathscr B_1,\mathscr B_2$ be translation invariant BF-spaces, $\fy \in C^\infty(\rr{d})$, $\omega _1,\omega _2\in
\mathscr{P}(\rr{2d})$ and $f\in \mathscr{D}'(X)$. If $\FB _1(\omega _1)_{loc}\subseteq \FB _2(\omega _2)_{loc}$, then
\begin{equation*}
\WF _{\FB _2(\omega _2)}(\varphi f)\subseteq \WF _{\FB _1(\omega
_1)}(f).
\end{equation*}
\end{thm}

\par

\begin{proof}
It suffices to prove
\begin{equation}\label{chi-subsetAA}
\Sigma_{ {\mathcal B _2} } (\fy f) \subseteq
\Sigma_{\mathcal B_1} (f).
\end{equation}
when $\mathcal B_j = \FB _j(\omega _j)$, $ \fy \in  \mathscr S(\rr d)$ and $f\in
\mathscr E'(\rr d)$, since the statement only involve local assertions. The local properties and Remark \ref{newbfspaces} also imply that it is no restriction to assume that $\omega _1 =\omega _2 = 1$.

\par

Let $ \xi_0 \in \Theta_{\mathcal B_2} (f)$,
and choose open cones $\Gamma _1$ and $\Gamma_2$ in $\rr d$ such that
$\overline {\Gamma _2}\subseteq \Gamma _1$. Since $f$ has compact
support, it follows that $|\widehat f(\xi )|\le
C\eabs \xi ^{N_0}$ for some positive constants $C$ and $N_0$. The result
therefore follows if we prove that for each $N$, there are constants
$C_N$ such that
\begin{multline}\label{cuttoff1}
|\fy f|_{ {\mathcal B _2} (\Gamma _2)}\le C_N \Big (|
f|_{\mathcal B_1(\Gamma _1)} +\sup _{\xi \in \rr
d} \big ( |\widehat f(\xi )|\eabs \xi ^{-N} \big )
\Big )
\\[1ex]
\text{when}\quad \overline \Gamma _2\subseteq \Gamma
_1\quad \text{and}\quad N=1,2,\dots .
\end{multline}

\par

By using the fact that $\omega$ is $v$-moderate for some $v\in
\mathscr P(\rr d) $, and letting $F(\xi )=|\widehat f(\xi )|$ and $\psi (\xi )=|\widehat \fy (\xi )|$, it
follows that $\psi$ turns rapidly to zero at infinity and
\begin{multline*}
|\fy f| _{ {\mathcal B _2}(\Gamma _2)} = |\varphi f|_{\FB _2(\Gamma_2)}
=
\|\mathscr{F}(\varphi f) \chi_{\Gamma _2}\|_{\mathscr{B}_2}
\\[1ex]
\leq
C \Big \| \Big ( \int_{\rr{d}} \widehat \fy (\cdo - \eta )\widehat
f(\eta )\, d\eta \Big ) \chi_{\Gamma_2} \Big \|_{\mathscr{B}_2}
\leq
C(J_1 + J_2)
\end{multline*}
for some positive constant $C$, where
\begin{align}
J_1
&=
\Big \| \Big ( \int _{\Gamma _1} \widehat \fy (\cdo - \eta)\widehat
f(\eta )\, d\eta \Big ) \chi _{\Gamma_2} \Big \| _{\mathscr{B}_2}\label{J1def}
\intertext{and}
J_2
&=
\Big \| \Big ( \int _{\complement \Gamma _1}\widehat \fy (\cdo -
\eta)\widehat f (\eta )\, d\eta \Big )\chi _{\Gamma _2} \Big
\| _{\mathscr{B} _2}\label{J2def}
\end{align}
and $\chi_{\Gamma_2}$ is the characteristic function of
$\Gamma_2$. First we estimate $J_1$. By (3) in Definition
\ref{BFspaces} and \eqref{propupps}, it follows for some constants
$C_1,\dots ,C_5$ that
\begin{multline}\label{J1comp}
J_1\leq  C_1 \Big \| \int _{\Gamma _1}\widehat \fy (\cdo -\eta )
\widehat f(\eta )\, d\eta \Big \| _{\mathscr{B}_2}
=
C_1 \nm {\widehat \fy * (\chi_{\Gamma_1}\widehat f) }{\mathscr{B}_2}
\\[1ex]
=
C_2 \nm {\fy  \mathscr F^{-1}(\chi_{\Gamma_1}\widehat f) }{\FB _2}
\le
C_3 \nm {\fy  \mathscr F^{-1}(\chi_{\Gamma_1}\widehat f) }{\FB _1}
\\[1ex]
=
C_4 \nm {\widehat \fy * (\chi_{\Gamma_1}\widehat f) }{\mathscr{B}_1}
\leq
C_5\nm {\widehat \fy \|_{L^1_{(v)}}\|\chi_{\Gamma_1}\widehat f
}{\mathscr{B}_1}
=
C_{\psi}|f|_{\FB _1(\Gamma_1)},
\end{multline}
where $C_{\psi}  = C_5\nm {\widehat \fy}{L^1_{(v)}}<\infty$, since
$\widehat \fy$ turns rapidly
to zero at infinity. In the second inequality we have used the fact
that $(\FB _1)_{loc}\subseteq (\FB _2)_{loc}$.

\par

In order to estimate $J_2$, we note that the conditions $\xi \in
\Gamma _2$, $\eta \notin \Gamma _1$ and the fact that $\overline
{\Gamma _2}\subseteq \Gamma _1$ imply that  $|\xi -\eta |>c\max
(|\xi|,|\eta |)$ for some constant $c>0$, since this is true when
$1=|\xi |\ge |\eta|$. We also note that if $N_1$ is large enough, then
$\eabs \cdo ^{-N_1}\in \mathscr B_2$, because $\mathscr S$ is
continuously embedded in $\mathscr B_2$. Since $\psi$ turns rapidly to
zero at infinity, it follows that for each $N_0> d+N_1$ and $N\in
\mathbf N$ such that $N > N_0$, it holds
\begin{multline}\label{J2comp}
J_2
\leq
C_1\Big \| \Big (\int_{\complement\Gamma _1}\eabs{\cdo -
\eta}^{-(2N_0+N)} F(\eta)\, d\eta \Big )\chi_{\Gamma_2}\Big \|
_{\mathscr{B}_2}
\\[1ex]
\leq
C_2\Big \| \Big (\int_{\complement \Gamma_1} \eabs \cdo^{-N_0} \eabs
\eta ^{-N_0}(\eabs \eta^{-N}F(\eta ))\, d\eta \Big
)\chi_{\Gamma_2}\Big \| _{\mathscr{B}_2}
\\[1ex]
\leq
C_2\int_{\complement \Gamma_1}\|\eabs \cdo
^{-N_0}\chi_{\Gamma_2}\|_{\mathscr{B}_2} \eabs \eta ^{-N_0}(|\eabs \eta
^{-N}F(\eta )|)\, d\eta
\\[1ex]
\leq
 C \sup_{\eta \in \rr{d}}|\eabs \eta^{-N}F(\eta )|,
\end{multline}
for some constants $C_1$ and $C_2 > 0$, where $C=C_2\nm {\eabs \cdo
^{-{N_0}}}{\mathscr B _2}\nm {\eabs \cdo ^{-{N_0}}}{L^1}<\infty$. This
proves \eqref{cuttoff1}, and the result follows.
\end{proof}

\par

\section{Mapping properties for pseudo-differential
operators on wave-front sets}\label{sec4}

\par

In this section we establish mapping properties for
pseudo-differential operators on wave-front sets of Fourier Banach
types. More precisely, we prove the following result (cf. \eqref{eq:inclusions}):

\par

\begin{thm}\label{mainthm2}
Let $\rho >0$, $\omega \in \mathscr P(\rr
{2d})$, $\omega _0 \in \mathscr P_{\rho ,0}(\rr
{2d})$, $a\in S^{(\omega _0)}_{\rho ,0} (\rr {2d})$, and $f\in \mathscr
S'(\rr d)$. Also let $\mathscr B$ be a translation
invariant BF-space on $\rr d$.
Then
\begin{multline}\label{wavefrontemb1}
\WF _{\FB (\omega /\omega _0)} (\op (a)f) \subseteq
\WF _{\FB (\omega )} (f)\\[1 ex]
\subseteq \WF _{\FB (\omega /\omega _0)} (\op (a)f)\ttbigcup
\Char _{(\omega _0 )}(a).
\end{multline}
\end{thm}

\par

We shall mainly follow the proof of Theorem 3.1 in \cite{PTT1}. The following restatement of Proposition 3.2 in \cite{PTT1} shows that  $ (x_0, \xi) \not\in
\WF _{\FB (\omega /\omega _0)} (\op (a)f) $ when $ x_0 \not\in \supp f $.

\par

\begin{prop}\label{propmain1AA}
Let $\omega \in \mathscr P(\rr {2d})$,  $\omega _0 \in \mathscr
P_{\rho ,\delta} (\rr {2d})$,
$0\le \delta \le \rho$, $0<\rho$, $\delta <1$, and let $a\in
S^{(\omega _0 )}_{\rho ,\delta}(\rr {2d})$. Also let $\mathscr B$ be a translation invariant BF-space, and let the operator
$L_a$ on $\mathscr S'(\rr d)$ be defined by the formula
\begin{equation}\label{Ladef}
(L_af)(x) \equiv  \fy _1(x)(\op (a)(\fy _2f))(x), \quad f\in \mathscr
S'(\rr d),
\end{equation}
where $\fy _1\in
C_0^{\infty}(\rr d)$ and $\fy _2 \in S_{0,0}^0(\rr d)$ are such that
$$
\supp \fy _1\bigcap \supp \fy _2=\emptyset .
$$
Then the kernel of $L_a$ belongs to
$\mathscr S(\rr {2d})$. In particular, the following is true:
\begin{enumerate}
\item $L_a =\op (a_0)$ for some $a_0\in \mathscr S(\rr {2d})$;

\vrum

\item $\WF _{\FB (\omega /\omega _0)}(L_af)=\emptyset$.
\end{enumerate}
\end{prop}

\par

Next we consider properties of the wave-front set of $\op (a)f$ at a
fixed point when $f$ is concentrated to that point.

\par

\begin{prop}\label{keyprop2AA}
Let $\rho$, $\omega$, $\omega _0$, $a$ and $\mathscr B$ be as in Theorem \ref{mainthm2}. Also let 
$f\in \mathscr E'(\rr d)$. Then the following is true:
\begin{enumerate}
\item if $\Gamma _1,\Gamma _2\subseteq \rr
d\back 0$ are open cones such that $\overline{\Gamma _2}\subseteq \Gamma _1$, and
$|f|_{\FB (\omega ,\Gamma _1)}<\infty$, then $|\op (a)f|_{\FB (\omega /\omega _0,\Gamma _2)}<\infty$;

\vrum

\item $\WF  _{\FB (\omega /\omega _0)}(\op (a)f)\subseteq
\WF _{\FB (\omega )}(f)$.
\end{enumerate}
\end{prop}

\par

We note that $\op (a)f$ in Proposition \ref{keyprop2AA} makes sense as
an element in $\mathscr S'(\rr d)$, by Proposition \ref{p5.4}.

\par

\begin{proof}
We shall mainly follow the proof of Proposition 3.3 in
\cite{PTT1}. We may assume that $\omega (x,\xi )=\omega (\xi )$, $\omega
_0(x,\xi )=\omega _0(\xi )$, and that $\supp a\subseteq K\times \rr d$
for some compact set $K\subseteq \rr d$, since the
statements only involve local assertions.

\par

Let  $F(\xi )=|\widehat f(\xi )\omega (\xi )|$, and let $\mathscr
F_1a$ denote the partial Fourier transform of $a(x, \xi)$ with respect
to the $x$ variable. By straightforward computation, for arbitrary $N$ we have
\begin{equation}\label{estpseudo}
|\mathscr{F}(\op(a)f)(\xi)\omega(\xi)/\omega _0(\xi)|
\leq
 C \int_{\rr{d}} \eabs {\xi-\eta}^{-N} F(\eta)\, d\eta ,
\end{equation}
for some constant $C$ (cf. (3.6) and (3.8) in \cite{PTT1}).

\par

We have to estimate
$$
|(\op(a)f)|_{\FB (\omega/\omega _0,\Gamma_2)}
=
\| \mathscr{F}(\op(a)f)\omega/\omega _0 \chi_{\Gamma_2}\| _{\mathscr{B}}.
$$
By \eqref{estpseudo} we get
\begin{multline*}
\|\mathscr{F}(\op(a)f)\omega/\omega _0 \chi_{\Gamma_2}\|_{\mathscr{B}}
\leq
C\Big \| \Big ( \int \eabs {\cdot - \eta}^{-N}F(\eta ) \, d\eta
\Big ) \chi _{\Gamma _2}\Big \| _{\mathscr{B}}
\cr
\leq
C(J_1 + J_2),
\end{multline*}
where $C$ is a constant and
\begin{align*}
J_1 &= \Big \| \Big ( \int _{\Gamma_1}\eabs {\cdot - \eta}^{-N}F(\eta
) \, d\eta \Big ) \chi _{\Gamma_2 }\Big \| _{\mathscr{B}}
\intertext{and}
J_2 &= \Big \| \Big ( \int _{\complement \Gamma_1}\eabs {\cdot -
\eta}^{-N}F(\eta ) \, d\eta \Big ) \chi _{\Gamma_2 }\Big \|
_{\mathscr{B}}.
\end{align*}

\par

In order to estimate $J_1$ and $J_2$ we argue as in the proof of
\eqref{cuttoff1}. More precisely, by \eqref{propupps} we get
\begin{multline*}
J_1 \leq
\Big \| \int_{\Gamma_1}\eabs {\cdo  - \eta}^{-N}F(\eta ) \, d\eta \Big
\|_{\mathscr{B}}
=
\|\eabs {\cdo }^{-N}* (\chi_{\Gamma_1}F)\| _{\mathscr{B}}
\\[1ex]
\leq
C\| \eabs \cdo ^{-N}\|_{L^1_{(v)}}\|\chi_{\Gamma_1}F\|_{\mathscr{B}}
<
\infty.
\end{multline*}

\par

Next we estimate $J_2$. Since $\overline {\Gamma_2} \subseteq \Gamma_1$, we get
$$
|\xi -\eta |\ge c\max ( |\xi |,|\eta |),\quad \text{when}\quad \xi \in
\Gamma _2,\ \text{and}\ \eta \in \complement \Gamma _1,
$$
for some constant $c>0$. (Cf. the proof of Proposition 3.3.)

\par

Since $f$ has compact support, it follows that $F(\eta )\le C\eabs
\eta ^{t_1}$ for some constant $C$. By combining these estimates we
obtain
\begin{multline*}
J_2
\leq
\Big \| \Big (\int _{\complement \Gamma_1}F(\eta)\eabs {\cdo -
\eta}^{-N} \, d\eta \Big ) \chi _{\Gamma _2}\Big \| _{\mathscr{B}}
\\[1ex]
\leq
C\Big \| \Big ( \int_{\complement \Gamma_1}\eabs {\eta}^{t_1}\eabs
\cdo ^{-N/2} \eabs \eta^{-N/2} \, d\eta \Big ) \chi_{\Gamma_2}\Big \|
_{\mathscr{B}}
\\[1ex]
\leq
C\|\eabs \cdo ^{-N/2}\chi_{\Gamma _2}\| _{\mathscr{B}} \int
_{\complement \Gamma_1} \eabs \eta^{-N/2+t_1} \, d \eta.
\end{multline*}
Hence, if we choose $N$ sufficiently large, it follows that the
right-hand side is finite. This proves (1).

\par

The assertion (2) follows immediately from (1) and the
definitions. The proof is complete.
\end{proof}

\par

\begin{proof}[Proof of Theorem \ref{mainthm2}]
By Proposition \ref{pseudocomp} it is no restriction to assume that
$t=0$.
We start to prove the first inclusion in \eqref{wavefrontemb1}.
Assume that $(x_0,\xi _0)\notin \WF _{\FB {(\omega)}}(f)$, let $\chi \in C_0^\infty (\rr d)$ be such that $\chi =1$
in a neighborhood of $x_0$, and set $\chi _1=1-\chi$ and $a_0(x,\xi
)=\chi (x)a(x,\xi )$. Then it follows from Proposition
\ref{propmain1AA} that
$$
(x_0,\xi _0)\notin \WF _{\FB {(\omega /\omega _0)}}(\op (a)(\chi
_1f)).
$$
Furthermore, by Proposition \ref{keyprop2AA} we get
\begin{equation*}
(x_0,\xi _0)\notin \WF _{\FB {(\omega /\omega _0
)}}(\op (a_0)(\chi f)),
\end{equation*}
which implies that
$$
(x_0,\xi _0)\notin \WF _{\FB {(\omega /\omega _0 )}}(\op (a)(\chi f)),
$$
since $\op (a)(\chi f)$ is equal to $\op (a_0)(\chi f)$ near
$x_0$. The result is now a consequence of the inclusion
\begin{multline*}
\WF _{\FB {(\omega /\omega _0)}}(\op (a)f)
\\[1ex]
\subseteq \WF _{\FB {(\omega /\omega _0 )}}(\op (a)(\chi f))\ttbigcup \WF _{\FB
{(\omega /\omega _0 )}}(\op (a)(\chi _1f)).
\end{multline*}

\par

It remains to prove the last inclusion in \eqref{wavefrontemb1}. By
Proposition \ref{propmain1AA} it follows that it is no restriction to
assume that $f$ has compact support. Assume that
$$
(x_0,\xi _0)\notin  \WF _{\FB (\omega /\omega _0)}(\op
(a)f)\ttbigcup \Char _{(\omega _0)}(a),
$$
and choose $b$, $c$ and $h$ as in Proposition \ref{psiecharequiv} (4). We
shall prove that $(x_0,\xi _0)\notin  \WF  _{\FB (\omega )}(f)$. Since
$$
f = \op (1-c)f + \op (b)\op (a)f-\op (h)f,
$$
the result follows if we prove
$$
(x_0,\xi _0)\notin \mathsf S_1\ttbigcup \mathsf S_2\ttbigcup \mathsf
S_3,
$$
where
\begin{align*}
\mathsf S_1 &= \WF _{\FB (\omega )}(\op (1-c)f),\quad
\mathsf S_2 = \WF _{\FB (\omega )}(\op (b)\op (a)f)
\\[1ex]
\text{and}\quad \mathsf S_3 &= \WF _{\FB (\omega )}(\op
(h)f).
\end{align*}

\par

We start to consider $\mathsf S_2$. By the first embedding in
\eqref{wavefrontemb1} it follows that
$$
\mathsf S_2 = \WF _{\FB (\omega )}(\op (b)\op (a)f)
\subseteq \WF  _{\FB (\omega /\omega _0)}(\op (a)f).
$$
Since we have assumed that $(x_0,\xi _0)\notin \WF _{\FB (\omega /\omega _0)}(\op
(a)f)$, it follows that $(x_0,\xi _0)\notin \mathsf S_2$.

\par

Next we consider $\mathsf S_3$. Since $h\in \mathscr S$, it follows that $\op (h)f\in \mathscr S$. Hence $S_3$ is empty.

\par

Finally we consider $\mathsf S_1$. By the assumptions it follows
that $c_0=1-c$ is zero in $\Gamma$, and by replacing $ \Gamma $ with a
smaller cone, if necessary, we may assume that $ c_0 = 0 $ in a
conical neighborhood of $ \Gamma$. Hence, if $\Gamma \equiv \Gamma _1$, $\Gamma _2$, $J_1$ and $J_2$ are the
same as in the
proof of Proposition \ref{keyprop2AA}, then it follows from that proof
and the fact that $c_0(x,\xi )\in S^0_{\rho ,0}$ is compactly
supported in the $x$-variable, that $J_1<+\infty$, and that for each $N\ge
0$, there are constants $C_N$ and $C_N'$ such that
\begin{multline}\label{estagain4}
| \op (c_0)f |_{\FB (\omega /\omega _0,\Gamma _2)} \le
C_N (J_1+J_2)
\\[1ex]
\le C_N' \Big (J_1+ \Big \Vert \int _{\complement
\Gamma _1}\eabs {\cdo} ^{-N}\eabs \eta ^{-N}\, d\eta \, \chi _{\Gamma _2} 
\Big \Vert _{\mathscr B} \Big ).
\end{multline}
By choosing $N$ large enough, it follows that
$$
|\op (c_0)f |_{\FB(\omega/\omega _0,\Gamma _2)}
< \infty .
$$
This proves that $(x_0,\xi _0)\notin \mathsf S_1$, and
the proof is complete.
\end{proof}

\par

\begin{rem}\label{rho0}
We note that the statements in Theorems \ref{mainthm2} are not true if $\omega _0=1$ and the assumption $\rho >0$ is replaced by
$\rho =0$. (Cf. Remark 3.7 in \cite{PTT1}.)
\end{rem}

\medspace

Next we apply Theorem \ref{mainthm2} on
operators which are elliptic with respect to $S^{(\omega _0)}_{\rho
,\delta}(\rr {2d})$, where $\omega _0\in \mathscr P_{\rho,\delta}(\rr {2d})$. More
precisely, assume that $0\le \delta <\rho \le1$ and $a\in S^{(\omega
_0)}_{\rho ,\delta}(\rr {2d})$. Then $a$ and
$\op (a)$ are called (locally) \emph{elliptic} with respect to
$S^{(\omega _0)}_{\rho ,\delta}(\rr {2d})$ or $\omega _0$,
if for each compact set $K\subseteq \rr d$, there are positive
constants $c$ and $R$ such that
$$
|a(x,\xi )| \ge c\omega _0(x,\xi ),\quad x\in K,\ |\xi |\ge R.
$$
Since $|a(x,\xi )|\le C\omega _0(x,\xi )$, it follows from the
definitions that for each multi-index $\alpha$, there are constants
$C_{\alpha,\beta}$ such that
\begin{equation*}
|\partial ^\alpha _x\partial ^\beta _\xi a(x,\xi )| \le C_{\alpha
,\beta}|a(x,\xi )|\eabs \xi ^{-\rho |\beta |+\delta |\alpha |},\quad
x\in K,\   |\xi |>R,
\end{equation*}
when $a$ is elliptic. (See e.{\,}g. \cite {Ho1,BBR}.) 

\par

It immediately follows from the
definitions that $\Char _{(\omega _0)}(a)=\emptyset$ when $a$ is
elliptic with respect to $\omega _0$. The following result is now an
immediate consequence of Theorem \ref{mainthm2}.

\par

\begin{thm}\label{hypoellthm}
Let $\omega \in \mathscr P(\rr
{2d})$, $\omega _0\in \mathscr P_{\rho ,0}(\rr {2d})$, $\rho >0$, and
let $a\in S^{(\omega _0)} _{\rho ,0} (\rr
{2d})$ be elliptic with respect to $\omega _0$. Also let $\mathscr B$
be a translation invariant BF-space.
If $f\in \mathscr S'(\rr d)$, then
$$
\WF _{\FB (\omega /\omega _0)} (\op (a)f)= \WF _{\FB (\omega )} (f).
$$
\end{thm}

\par

\section{Wave-front sets of sup and inf types and
pseudo-differential operators}\label{sec5}

\par

In this section we put the micro-local analysis in a more general
context compared to previous sections,
and define wave-front sets with respect to sequences of
Fourier BF-spaces.

\par

Let $\omega _j \in \mathscr
P(\rr {2d})$ and $\mathscr B_j$ be translation invariant BF-space on
$\rr d$ when $j$ belongs to some index set $J$, and consider the array
of spaces, given by
\begin{equation}\label{notconvsequences}
(\mathcal B_j) \equiv (\mathcal B_j)_{j\in J},\quad
\text{where}\quad \mathcal B_j=\FB _j {(\omega
_j)}, \quad j \in J.
\end{equation}

\par

If $f\in \mathscr S'(\rr d)$, and $(\mathcal B_j)$ is given by
\eqref{notconvsequences}, then we let
$\Theta_{(\mathcal B_j) }^{\sup}(f)$ be the set of all $\xi \in \rr
d\back 0$ such that for some $\Gamma = \Gamma _{\xi}$ and \emph{each}
$j\in J$ it holds $|f|_{\mathcal B_j(\Gamma )} < \infty$. We
also let $\Theta _{(\mathcal B_j) }^{\inf}(f)$ be the set of all $
\xi \in \rr d\back 0 $ such that for some $\Gamma = \Gamma
_{\xi}$ and \emph{some} $j\in J$ it holds $|f|_{\mathcal B_j(\Gamma)}
<\infty$. Finally we let $\Sigma _{(\mathcal B_j) }^{\sup} (f)$ and
$\Sigma _{(\mathcal B_j) }^{\inf} (f)$ be the complements in $\rr
d\back 0 $ of $\Theta_{(\mathcal B_j) }^{\sup}(f)$ and $\Theta
_{(\mathcal B_j) }^{\inf} (f)$ respectively.

\par

\begin{defn}\label{defsuperposWF}
Let $J$ be an index set, $\mathscr B_j$ be translation invariant
BF-space on $\rr d$, $\omega _j\in \mathscr P(\rr {2d})$ when $j\in
J$, $(\mathcal B_j)$ be as in
\eqref{notconvsequences}, and let $X$ be an open subset of $\rr d$.
\begin{enumerate}
\item The wave-front set of $f\in \mathscr D'(X)$,
$\WF ^{\, \sup} _{(\mathcal B_j) }(f) = \WF ^{\, \sup} _{(\FB
_j(\omega _j) )}(f)$,
of \emph{sup-type} with respect to $(\mathcal B_j) $,
consists of all pairs
$(x_0,\xi_0)$ in $X\times (\rr d \back 0)$ such that
$
\xi _0 \in  \Sigma ^{\sup} _{(\mathcal B_j)} (\fy f)
$
holds for each $\fy \in C_0^\infty (X)$ such that $\fy (x_0)\neq
0$;

\vrum

\item The wave-front set of $f\in \mathscr D'(X)$,
$\WF ^{\, \inf} _{(\mathcal B_j) }(f) = \WF ^{\, \inf} _{(\FB
_j(\omega _j) )}(f)$, of \emph{inf-type} with respect to $(\mathcal
B_j)$, consists of all pairs
$(x_0,\xi_0)$ in $X\times (\rr d \back 0)$ such that
$
\xi _0 \in  \Sigma ^{\inf} _{(\mathcal B_j)} (\fy f)
$
holds for each $\fy \in C_0^\infty (X)$ such that $\fy (x_0)\neq
0$.
\end{enumerate}
\end{defn}

\par

\begin{rem}\label{remstandWF}
Let $\omega _j(x,\xi ) = \eabs \xi ^{-j}$ for $j\in J=\mathbf
N_0$ and $\mathscr B_j=L^{q_j}$, where $q_j\in [1,\infty ]$. Then it
follows that $\WF _{(\mathcal B_j)}^{\, \sup}(f)$ in
Definition \ref{defsuperposWF} is equal to the standard wave front
set $\WF (f)$ in Chapter VIII in \cite{Ho1}.
\end{rem}

\par

The following result follows immediately from Theorem \ref{mainthm2}
and its proof. We omit the details.

\par

\renewcommand{\rubrik}{Theorem \ref{mainthm2}$'$}

\begin{tom}
Let $\rho >0$, $\omega _j\in \mathscr P(\rr
{2d})$ for $j\in J$, $\omega _0 \in \mathscr P_{\rho ,0}(\rr
{2d})$, $a\in S^{(\omega _0)}_{\rho ,0} (\rr {2d})$ and $f\in \mathscr
S'(\rr d)$. Also let $\mathscr B_j$ be a translation
invariant BF-space on $\rr d$ for every $j$. Then
\begin{multline}\tag*{(\ref{wavefrontemb1})$'$}
\WF ^{\, \sup}_{(\FB _j(\omega _j/\omega _0))} (\op (a)f) \subseteq
\WF ^{\, \sup}_{(\FB_j (\omega _j))} (f)
\\[1ex]
\subseteq \WF ^{\, \sup}_{(\FB _j(\omega _j/\omega _0))} (\op
(a)f)\ttbigcup \Char _{(\omega _0)}(a),
\end{multline}
and
\begin{multline}\tag*{(\ref{wavefrontemb1})$''$}
\WF ^{\, \inf}_{(\FB _j(\omega _j/\omega _0))} (\op (a)f) \subseteq
\WF ^{\, \inf}_{(\FB_j (\omega _j))} (f)
\\[1ex]
\subseteq \WF ^{\, \inf}_{(\FB _j(\omega _j/\omega _0))} (\op (a)f)\ttbigcup
\Char _{(\omega _0)}(a).
\end{multline}
\end{tom}

\par

The following generalization of Theorem \ref{hypoellthm} is an
immediate consequence of Theorem \ref{mainthm2}$'$, since $\Char
_{(\omega _0)}(a)=\emptyset$, when $a$ is elliptic with respect to
$\omega_0$.

\par

\renewcommand{\rubrik}{Theorem \ref{hypoellthm}$'$}

\begin{tom}

Let $\rho >0$, $\omega _j\in \mathscr P(\rr
{2d})$ for $j\in J$, $\omega _0 \in \mathscr P_{\rho ,0}(\rr
{2d})$ and let $a\in S^{(\omega _0)}_{\rho ,0} (\rr {2d})$ be elliptic
with respect to $\omega _0$.  Also let $\mathscr B_j$ be a translation
invariant BF-space on $\rr d$ for every $j$. If $f\in \mathscr S'(\rr
d)$, then
$$
\WF ^{\, \sup}_{(\FB _j(\omega _j/\omega _0))} (\op (a)f) =
\WF ^{\, \sup}_{(\FB_j (\omega _j))} (f)
$$
and 
$$
\WF ^{\, \inf}_{(\FB _j(\omega _j/\omega _0))} (\op (a)f) =
\WF ^{\, \inf}_{(\FB_j (\omega _j))} (f).
$$

\end{tom}

\par

\begin{rem}\label{remhormWFsets}
We note that many properties valid for the wave-front sets of Fourier
Banach type also hold for wave-front sets in the present
section. For example, the conclusions in Remark \ref{rho0}
hold for wave-front sets of sup- and inf-types.
\end{rem}

\par

Finally we remark that there are some technical generalizations of
Theorem \ref{mainthm2} which involve pseudo-differential operators
with symbols in $S^{(\omega _0)}_{\rho ,\delta} (\rr {2d})$ with $0\le
\delta <\rho \le 1$. From these generalizations it follows that
$$
\WF (\op (a)f) \subseteq \WF(f) \subseteq \WF (\op (a)f)\ttbigcup \Char
_{(\omega _0)}(a),
$$
when $0\le \delta <\rho \le 1$, $\omega _0\in \mathscr P_{\rho ,\delta
}(\rr {2d})$, $a\in S^{(\omega _0)}_{\rho ,\delta}(\rr {2d})$ and
$f\in \mathscr S'(\rr d)$. (Cf. Theorem 5.3$'$ and Theorem 5.5 in
\cite{PTT1}.)

\par

\section{Wave front sets with respect to modulation spaces}\label{sec6}

\par

In this section we define wave-front sets with respect to modulation
spaces, and show that they coincide with wave-front sets of Fourier
Banach types. In particular, all micro-local properties for
pseudo-differential operators in the previous sections carry over to
wave-front sets of modulation space types.

\par

We start with defining general types of modulation spaces. Let
(the window) $\phi \in \mathscr S'(\rr d)\back  0$ be
fixed, and let $f\in \mathscr S'(\rr d)$. Then the short-time Fourier
transform $V_\phi f$ is the element in $\mathscr S'(\rr {2d})$,
defined by the formula
$$ 
(V_{\phi} f)(x,\xi) \equiv \mathscr{F}(f\cdot
\overline{\phi(\cdot-x)})(\xi).
$$
We usually assume that $\phi \in \mathscr{S}(\rr{d})$, and in
this case the short-time Fourier transform $(V_{\phi}f)$ takes
the form 
$$ 
(V_{\phi} f)(x,\xi)
=
(2\pi)^{-d/2}\int_{\rr{d}} f(y)\overline{\phi(y-x)}e^{-i\scal y \xi
}\, dy,
$$ 
when $f\in \mathscr{S}(\rr{d})$.

\par

Now let $\mathscr{B}$ be a translation invariant BF-space on
$\rr {2d}$, with respect to $v\in \mathscr P(\rr {2d})$. Also let
$\phi \in \mathscr{S}(\rr{d})\back{0}$ and $\omega \in
\mathscr{P}(\rr{2d})$ be such that $\omega$ is $v$-moderate. Then the
modulation space $M(\omega)=M(\omega ,\mathscr B)$ consists of all $f\in
\mathscr{S}'(\rr{d})$ such that $V_{\phi}f\cdot \omega \in
\mathscr{B}$. We note that $M(\omega ,\mathscr B)$ is a Banach space
with the norm
\begin{equation}\label{modnorm}
\|f\|_{M(\omega ,\mathscr B)} \equiv \|V_{\phi} f
\omega\|_{\mathscr{B}}
\end{equation}
(cf. \cite{Feichtinger3}).

\par

\begin{rem}\label{Modamalgam}
Assume that $p,q\in [1,\infty]$, $\omega\in \mathscr P (\rr {2d})$ and
let $L^{p,q}_{1}(\rr{2d})$ and $L^{p,q}_{2}(\rr{2d})$ be the sets of
all $F\in  L^1_{loc} (\rr{2d})$ such that
\begin{equation*}
\|F\|_{L^{p,q}_1}  \equiv \Big ( \int \Big( \int |F(x,\xi)|^p\,
dx\Big )^{q/p}\,d\xi \Big )^{1/q}
<\infty
\end{equation*}
and
\begin{equation*}
\|F\|_{L^{p,q}_2} \equiv \Big ( \int \Big ( \int |F(x,\xi)|^p\,
d\xi \Big )^{q/p}\, dx\Big )^{1/q}<\infty ,
\end{equation*}
respectively (with obvious modifications when $p=\infty$ or
$q=\infty$). Then $M(\omega ,\mathscr B)$ is equal to the usual
modulation space $M^{p,q}_{(\omega)}(\rr{d})$ when
$\mathscr{B}=L^{p,q}_1(\rr{2d})$. If instead
$\mathscr{B}=L^{p,q}_2(\rr{2d})$, then $M(\omega ,\mathscr B)$ is
equal to the space $W^{p,q}_{(\omega)}(\rr{d})$, related to
Wiener-amalgam spaces.
\end{rem}

\par

In the following proposition we list some important properties for
modulation spaces. We refer to \cite{Gro-book} for the proof.

\par

\begin{prop}\label{modproperties}
Assume that $\mathscr B$ is a translation invariant BF-space on $\rr
{2d} $with respect to $v\in \mathscr P(\rr {2d})$, and that
$\omega _0, v_0\in\mathscr{P}(\rr{2d})$ are such that $\omega$
is $v$-moderate. Then the following is true:
\begin{enumerate}
\item if $\phi\in M^1_{(v_0v)}(\rr{d})\back{0}$, then $f\in M(\omega
,\mathscr B)$ if and only if $V_\phi f \omega \in \mathscr
B$. Furthermore, \eqref{modnorm} defines a norm on $M(\omega
,\mathscr B)$, and different choices of $\phi$ gives rise to equivalent
norms;

\vrum

\item $ M^1_{(v_0v)}(\rr d)\subseteq M(\omega ,\mathscr{B})\subseteq
M^{\infty}_{(1/(v_0v))}(\rr d)$.
\end{enumerate}
\end{prop}

\par

The following generalization of Theorem 2.1 in \cite{RSTT} shows that
modulation spaces are locally the same as translation invariant
Fourier BF-spaces. We recall that if $\fy \in \mathscr S  (\rr d)\back
0$ and $\mathscr B$ is a translation invariant BF-space  on $\rr
{2d}$, then it follows from Proposition \ref{propbnoll} that
\begin{equation}\label{B0def}
\mathscr{B}_0 \equiv \sets {f\in \mathscr
S'(\rr d)}{\fy \otimes f \in \mathscr B}
\end{equation}
is a translation invariant BF-space on $\rr d$ which is independent of
the choice of $\fy$.

\par

\begin{prop}\label{propekvnorm}
Let $\fy \in C_0^\infty (\rr d)\back 0$,
$\mathscr B$ be a translation invariant BF-space  on $\rr {2d}$, and
let $\mathscr B_0$ be as in \eqref{B0def}. Also let $\omega \in
\mathscr P(\rr {2d})$, and $\omega _0(\xi )=\omega (x_0,\xi )$,
for some fixed $x_0\in \rr d$. Then
$$
M(\omega ,\mathscr B)\cap \mathscr E'(\rr d) = \FB _0(\omega _0)\cap
\mathscr E'(\rr d).
$$
Furthermore, if $K\subseteq \rr d$ is compact, then 
\begin{equation}\label{locest1}
C^{-1}\nm {f}{\FB _0(\omega _0)}\le \nm f{M(\omega ,\mathscr B)}\le
C\nm {f}{\FB _0(\omega _0)},\quad f\in \mathscr E'(K),
\end{equation}
for some constant $C$ which only depends on $K$.
\end{prop}

\par

We need the following lemma for the proof.

\par

\begin{lemma}\label{korttidslemma}
Assume that $f\in \mathscr{E}'(\rr{d})$. Then the following is true:
\begin{enumerate}
\item \label{kortl1} if $\phi\in C^{\infty}_0(\rr{d})$, then there
exists $0\le \fy \in C^{\infty}_0(\rr{d})$ such that
\begin{equation}\label{stft-ft}
(V_{\phi}f)(x,\xi )=\fy (x)(\widehat {f}*(\mathscr {F}(\overline
{\phi(\cdot-x)})))(\xi )\, \text ;
\end{equation}

\vrum

\item \label{kortl2} if $\fy \in C^{\infty}_0(\rr{d})$, then there
exists $\phi \in C^{\infty}_0(\rr{d})$ such that
\begin{equation}\label{ft-stft}
(\fy \otimes \widehat{f}) (x,\xi) =\fy (x)V_{\phi}f(x,\xi ).
\end{equation}
\end{enumerate}
\end{lemma}

\par

\begin{proof}
(1) Let $\fy \in C_0^\infty$ be equal to $(2\pi )^{d/2}$ in a
compact set containing the support of the map $x\mapsto V_\phi f(x,\xi
)$. Then (1) is a straight-forward consequence of Fourier's
inversion formula.

\par

The assertion (2) follows by choosing $\phi \in C_0^\infty$ such that
$\phi =1$ on $\supp f -\supp \fy$. 
\end{proof}

\par

\begin{proof}[Proof of Proposition \ref{propekvnorm}]
We may assume that $\omega =\omega _0=1$ in view of Remark
\ref{newbfspaces}. Assume that $f\in
\mathscr E'$ and $\fy \in C_0^\infty (\rr d)\back 0$. From
\eqref{kortl2} of Lemma \ref{korttidslemma} it follows that there
exists $\phi \in C^{\infty}_0$ such that
\begin{equation*}
\| f\|_{M(\mathscr{B})}= \| V_{\phi}f \| _{\mathscr{B}}=
\|\varphi \otimes \widehat{f}\| _{\mathscr{B}}= \nm {\widehat f}{\mathscr
B_0},
\end{equation*}
and \eqref{locest1} follows. The proof is complete.
\end{proof}

\medspace

Let $\mathscr B$ be a translation invariant BF-space on $\rr
{2d}$, $\phi \in \mathscr{S}(\rr{d} )\back{0}$ be fixed,
$\omega \in \mathscr{P}(\rr{2d})$, $\Gamma \subseteq \rr{d}\back{0}$
be an open cone, and let $\chi_{\Gamma}(x,\xi)=\chi_{\Gamma}(\xi)$ be
the characteristic function of $\Gamma$. For any $f\in
\mathscr{S}'(\rr{d})$ we set
\begin{multline}\label{modseminorm}
|f|_{\mathcal B(\Gamma )} =
|f|_{\modBrk{}{\Gamma}{}}  =|f|_{\modBrk{\phi}{\Gamma}{}}
\equiv
\|(V_{\phi}f) \omega\chi_{\Gamma}\|_{\mathscr{B}}
\\[1ex]
\text{when}\quad \mathcal B=M{(\omega ,\mathscr B)}.
\end{multline}
We note that $|\cdo |_{\mathcal B(\Gamma )}$ defines a semi-norm
on $\mathscr S'$ which might attain the value $+\infty$. If $\Gamma
=\rr d\back 0$, then $|f|_{\mathcal B(\Gamma )} = \nm
f{M{(\omega ,\mathscr B)}}$.

\par

Let $\mathscr B$ be a translation invariant BF-space on $\rr {2d}$,  $\omega \in \mathscr P(\rr {2d})$, $f\in \mathscr D'(X)$, and let
$\mathcal B=M{(\omega ,\mathscr B)}$. Then $\Theta _{\mathcal B}(f)$, $\Sigma
_{\mathcal B}(f)$ and the wave-front set $\WF
_{\mathcal B}(f)$ of $f$ with respect to the modulation space
$\mathcal B$ are defined in the same way as in Section
\ref{sec3}, after replacing the semi-norms of Fourier Banach types in
\eqref{notconv} with the semi-norms in \eqref{modseminorm}.

\par

In Theorem \ref{WFidentity} below we prove that wave-front sets of Fourier BF-spaces
and modulation space types agree with each others. As a
first step we prove that
$\WF_{\modBr{\phi}{}{}}(f)$ is independent of $\phi$ in \eqref{modseminorm}.

\par

\begin{prop}\label{fonsteroberoende}
Let $X\subseteq \rr d$ be open,
$f\in \mathscr{D}'(X)$ and $\omega\in \mathscr{P}(\rr{2d})$. Then $\WF_{\modBr{\phi}{}{}}(f)$ is
independent of the window function $\phi \in \mathscr S (\rr{d})
\back{0}$.
\end{prop}

\par

We need some preparation for the proof, and start with the following
lemma. We omit the proof (the result can be found in \cite{CG1}).

\par

\begin{lemma}\label{STFTdecay}
Let $f\in \mathscr E'(\rr d)$ and $\phi \in \mathscr S (\rr d)$. Then
for some constant $N_0$ and every $N\geq 0$, there are constants $C_N$
such that
$$
|V_{\phi} f(x,\xi)| \leq C_N \eabs x ^{-N}\eabs \xi^{N_0}.
$$
\end{lemma} 

\par

The following result can be found in \cite{Gro-book}. Here $\widehat
*$ is the twisted convolution, given by the formula
$$
(F\, \widehat *\, G)(x,\xi )=(2\pi )^{-d/2}\iint F(x-y,\xi -\eta
)G(y,\eta )e^{-i\scal {x-y}\eta}\, dyd\eta ,
$$
when $F,G\in \mathscr S(\rr {2d})$. The definition of $\widehat *$
extends in such way that one may permit one  of $F$ and $G$ to belong
to $\mathscr S'(\rr {2d})$, and in this case it follows that $F \,
\widehat * \, G$ belongs to $\mathscr S'\cap C^\infty$.

\par

\begin{lemma}\label{stftproperties}
Let $f\in \mathscr S'(\rr d)$ and $\phi _j\in
\mathscr S(\rr d)$ for $j=1,2,3$. Then
$$
(V_{\phi _1}f)\widehat*(V_{\phi _2}\phi _3) = (\phi
_3,\phi _1)_{L^2}\cdot V_{\phi _2}f.
$$
\end{lemma}

\par

\begin{proof}[Proof of Proposition \ref{fonsteroberoende}]
We assume that $f\in \mathscr{E}'(\rr d)$ and that $\omega(x,\xi) =\omega(\xi)$, since the statements only involve local assertions. Assume that $\phi, \phi_1 \in \mathscr S(\rr d)\back 0$ and let $\Gamma _1$ and $\Gamma_2$ be open cones in $\rr d$ such that $\overline{\Gamma_2}\subseteq \Gamma_1$. The assertion follows if we prove that 
\begin{equation}\label{modseminormineq}
|f|_{\modBrk{\phi}{\Gamma_2}{}} \leq C(|f|_{\modBrk{\phi_1}{\Gamma_1}{}} + 1)
\end{equation}
for some constant $C$.

\par

When proving \eqref{modseminormineq} we shall mainly follow the proof
of \eqref{cuttoff1}. Let
$v \in \mathscr P$ be chosen such that $\omega$ is $v$-moderate, and let
$$
\Omega_1 = \{ (x, \xi) ; \xi\in \Gamma_1 \} \subseteq \rr {2d} \qquad
\text{and} \qquad \Omega_2 = \complement \Omega_1 \subseteq \rr {2d},
$$
with characteristic functions $\chi_1$ and $\chi_2$ respectively. Also
set
$$
F _k(x, \xi) = |V_{\phi_1}f(x, \xi)\omega(\xi)\chi _k(x, \xi)| \qquad
\text{and} \qquad G = |V_{\phi}\phi_1(x, \xi)v(\xi)|.
$$
By Lemma \ref{stftproperties}, and the fact that $\omega$ is $v$-moderate we get
$$
|V_{\phi}f(x, \xi)\omega(x, \xi)| \leq  C((F_1 + F_2) * G)(x, \xi),
$$
for some constant $C$, which implies that
\begin{equation}\label{uppskattning}
|f|_{\modBrk{\phi}{\Gamma_2}{}} \leq C (J _1 + J _2), 
\end{equation}
where
$$
J_k = \nm{(F_k * G)\chi_{\Gamma_2}}{\mathscr B}
$$
and $\chi_{\Gamma_2}(x,\xi)=\chi_{\Gamma_2}(\xi)$ is the characteristic function of $\Gamma_2$.
Since $G$ turns rapidly to zero at infinity, \eqref{propupps} gives
\begin{equation}\label{del1}
J_1 \leq \nm{F_1*G}{\mathscr B} \leq \nm G{L^1_{(v)}} \nm {F_1}{\mathscr B}=C|f|_{\modBrk{\phi_1}{\Gamma_1}{}}, 
\end{equation}
where $C = \nm{G}{L^1_{(v)}}$.

\par

Next we consider $J_2$. Since, for each $N \geq 0$, there are constants $C_N$ such that
$$
F_2(x, \xi)=0,\qquad \text{and} \qquad \eabs{\xi-\eta}^{-2N}\leq C_N \eabs \xi^{-N}\eabs \eta^{-N}
$$
when $\xi\in \Gamma_2$ and $\eta \in \complement \Gamma_1$, it follows from Lemma \ref{STFTdecay} and the computations in \eqref{J2comp} that
$$
(F_2 * G)(x, \xi) \leq C_N \eabs x^{-N}\eabs \xi^{-N},\qquad  \xi\in \Gamma_2.
$$
Consequently, $J_2 < \infty$. The estimate \eqref{modseminormineq} is now a consequence of
\eqref{uppskattning}, \eqref{del1} and the fact that $J_2 < \infty$.
This completes the proof.
\end{proof}

\par

Since $\WF_{M^\phi(\omega,\mathscr B)} (f)$ is independent of $\phi$ we usually omit $\phi$ and write $\WF_{M(\omega,\mathscr B)} (f)$ instead. We are now able to prove the following.

\par

\begin{prop}\label{theta-sigma}
Assume that $\mathscr B$ is a translation invariant BF-space on $\rr {2d}$, $\mathscr B_0$ is given by \eqref{B0def}, $\phi\in \mathscr S (\rr d) \back{0}$ and
$\omega \in \mathscr P(\rr{2d})$. Also assume that $f\in
\mathscr E'(\rr{ d})$. Then
\begin{equation}\label{fy-invariant2}
\Theta_{\modBr{\phi} {}{}} (f) = \Theta_{\mathscr{FB}_0
(\omega )} (f)\quad \text{and}\quad
\Sigma_{\modBr{\phi} {}{}}  (f) = \Sigma_{\mathscr{FB}_0(\omega )} (f).
\end{equation}
\end{prop}

\par

\begin{proof}
We may assume that $\omega =1$ in view of Lemma \ref{newbfspaces}. Let $\Gamma _1,\Gamma _2$ be open cones in $\rr d\back 0$ such that $\overline {\Gamma _2}\subseteq \Gamma _1$, let $\chi _{\Gamma _2}(x,\xi )=\chi _{\Gamma _2}(\xi )$ be the characteristic function of $\Gamma _2$, and let $\fy$ and $\phi$ be chosen such that (1) in Lemma \ref{korttidslemma} is fulfilled.

\par

By \eqref{stft-ft} it follows that
$$
|V_\phi f(x,\xi )| \le \fy (x)(|\widehat f|*|\mathscr F\check \phi |)(\xi ).
$$
This gives
\begin{multline*}
|f|_{M^\phi (\Gamma _2,\mathscr B)} = |V_\phi f\chi _{\Gamma _2}|_{\mathscr B}
\le C | \fy \otimes \big (  (|\widehat f|*|\mathscr F\check \phi |) \chi _{\Gamma _2}\big )|_{\mathscr B}
\\[1ex]
= C | (|\widehat f|*|\mathscr F\check \phi |) \chi _{\Gamma _2} |_{\mathscr B_0}
\le C(J_1+J_2),
\end{multline*}
for some constant $C$, where $J_1$ and $J_2$ are the same as in \eqref{J1def} and \eqref{J2def} with $\mathscr B_2=\mathscr B_0$, $\psi =|\mathscr F\check \phi |$ and $F=|\widehat f|$.

\par

A combination of the latter estimate, \eqref{J1comp} and \eqref{J2comp} now gives that for each $N\ge 0$, there is a constant $C_N$ such that
$$
|f|_{M^\phi (\Gamma _2,\mathscr B)}  \le C_N\Big ( |f|_{\FB_0} +\sup _\xi |\widehat f(\xi )\eabs \xi ^{-N}|\Big ).
$$
Hence, by choosing $N$ large enough it follows that $|f|_{M^\phi (\Gamma _2,\mathscr B)}$ is finite when $ |f|_{\FB_0}<\infty$. Consequently,
\begin{equation}\label{thetaFBMB}
\Theta_{\mathscr{FB}_0} (f)\subseteq \Theta_{M(\mathscr B)} (f).
\end{equation}

\par

In order to get a reversed inclusion we choose $\fy$ and $\phi$ such that Lemma \ref{korttidslemma} (2) is fulfilled. Then \eqref{ft-stft} gives
\begin{multline*}
|f|_{\FB _0(\Gamma )} = \| \fy \otimes (\widehat{f} \, \chi _\Gamma ) \| _{\mathscr B}
=\| (\fy \otimes 1)(V_{\phi}f \, \chi _\Gamma ) \| _{\mathscr B} 
\\[1ex]
\le C_1\nm {\fy}{L^\infty} \nm {V_{\phi}f \, \chi _\Gamma}{\mathscr B}
= C_2 |f|_{M^\phi (\Gamma ,\mathscr B)},
\end{multline*}
for some constants $C_1, C_2 >0.$ This proves that \eqref{thetaFBMB} holds with reversed inclusion. The proof is complete.
\end{proof}

\par

The following result is now an immediate consequence of Proposition \ref{theta-sigma}.

\par

\begin{thm}\label{WFidentity}
Assume that $\mathscr B$ is a translation invariant BF-space on $\rr {2d}$, $\mathscr B_0$ is given by \eqref{B0def}, $\omega \in \mathscr P(\rr {2d})$, $X\subseteq \rr d$ is open and that $f\in \mathscr D'(X)$. Then
$$
\WF_{\FB _0(\omega )}(f) = \WF _{M(\omega ,\mathscr B)}(f).
$$
\end{thm}

\par


\vspace{2cm}

\begin{thebibliography}{150}

\bibitem{BaC}
W. Baoxiang, H. Chunyan \emph{Frequency-uniform decomposition
method for the generalized BO, KdV and NLS equations}, {J.
Differential Equations}, \textbf{239} (2007), 213--250.

\bibitem{BBR} P. Boggiatto, E. Buzano, L. Rodino \emph{Global
Hypoellipticity and Spectral Theory},  Mathematical Research, 92,
Akademie Verlag, Berlin, 1996.

\bibitem{CG1} E.~{C}ordero, K.~{G}r{\"o}chenig,
\emph{{T}ime-frequency analysis of localization operators},
 {J}. {F}unct. {A}nal., \textbf{205(1)} (2003), 107--131.

\bibitem{Czaja} W. Czaja, Z. Rzeszotnik
\emph{Pseudodifferential operators and Gabor frames: spectral
asymptotics}, {Math. Nachr.} \textbf{233-234} (2002), 77--88.

\bibitem{F1}  H.~G.~Feichtinger \emph{Modulation spaces on locally
compact abelian groups. Technical report}, {University of
Vienna}, Vienna, 1983; also in: M. Krishna, R. Radha,
S. Thangavelu (Eds) Wavelets and their applications, Allied
Publishers Private Limited, NewDehli Mumbai Kolkata Chennai Hagpur
Ahmedabad Bangalore Hyderbad Lucknow, 2003, pp. 99--140.

\bibitem{F2} \bysame \emph{Wiener amalgams over Euclidean spaces and some of their applications},
in: Function spaces (Edwardsville, IL, 1990), Lect. Notes in pure and
appl. math., 136, Marcel Dekker, New York, 1992, pp. 123–137.

\bibitem{Feichtinger3}  {H. G. Feichtinger and K. H. Gr{\"o}chenig}
\emph{Banach spaces related to integrable group representations and
their atomic decompositions, I}, J. Funct. Anal., \textbf{86}
(1989), 307--340.

\bibitem{Feichtinger4} \bysame \emph{Banach spaces related to
integrable group representations and their atomic decompositions, II},
Monatsh. Math., \textbf{108} (1989), 129--148.

\bibitem{Feichtinger5} \bysame \emph{Gabor frames and time-frequency
analysis of distributions}, {J. Functional
Anal.,} \textbf {146} (1997), 464--495.

\bibitem{Feichtinger6} \bysame \emph{Modulation spaces: Looking back and ahead},
Sampl. Theory Signal Image Process. \textbf{5} (2006), 109--140.

\bibitem{Fo}  {G. B. Folland} \emph
{Harmonic analysis in phase space}, {Princeton U. P., Princeton},
1989.

\bibitem{Grobner} P. Gr{\"o}bner \emph{Banachr{\"a}ume Glatter
Funktionen und Zerlegungsmethoden}, Thesis, University of Vienna,
Vienna, 1992.

\bibitem{Grochenig0a} {K. H. Gr{\"o}chenig} \emph {Describing
functions: atomic decompositions versus frames},
{Monatsh. Math.},\textbf{112} (1991), 1--42.

\bibitem{Gro-book} K. Gr\"{o}chenig, \newblock \textit{Foundations of
Time-Frequency Analysis},
\newblock Birkh\"auser, Boston, 2001.

\bibitem{Grochenig2} \bysame \emph{Composition and spectral invariance
of pseudodifferential operators on modulation spaces}, J. Anal.
Math., \textbf{98} (2006), 65--82.

\bibitem{Grochenig0}  {K. H. Gr{\"o}chenig and C. Heil} \emph
{Modulation spaces and pseudo-differential operators}, Integral
Equations Operator Theory, \textbf{34} (1999), 439--457.

\bibitem{Grochenig1b}  \bysame \emph {Modulation spaces as symbol
classes for pseudodifferential operators {\rm {in: M. Krishna,
R. Radha, S. Thangavelu (Eds) Wavelets and their applications}}},
Allied Publishers Private Limited, NewDehli Mumbai Kolkata Chennai
Hagpur Ahmedabad Bangalore Hyderbad Lucknow, 2003, pp. 151--170.

\bibitem{Grochenig1c} \bysame \emph{Counterexamples for boundedness of
pseudodifferential operators}, Osaka  J. Math., \textbf{41} (2004),
681--691.

\bibitem{Grochenig2a} K. Gr{\"o}chenig, M. Leinert \emph{Wiener's lemma
for twisted convolution and Gabor frames}, J. Amer. Math. Soc.,
\textbf{17} (2004), 1--18.

\bibitem{Herau1} F. H{\' e}rau \emph{Melin--H{\"o}rmander inequality
in a Wiener type pseudo-differential algebra}, Ark. Mat., \textbf{39}
(2001), 311--38.

\bibitem{HTW} A. Holst, J. Toft, P. Wahlberg \emph{Weyl product
algebras and modulation spaces}, {J. Funct. Anal.}, \textbf{251}
(2007), 463--491.

\bibitem{Ho1}  L. H{\"o}rmander \emph{The Analysis of Linear
Partial Differential Operators}, vol {I--III},
Springer-Verlag, Berlin Heidelberg NewYork Tokyo, 1983, 1985.

\bibitem{Hrm-nonlin} \bysame \emph{Lectures on Nonlinear Hyperbolic
Differential Equations}, Springer-Verlag, Berlin, 1997.

\bibitem{Okoudjou} K. Okoudjou \emph{Embeddings of some classical
Banach spaces into modulation spaces}, {Proc. Amer. Math. Soc.},
\textbf{132} (2004), 1639--1647.

\bibitem{Pilipovic2} {S. Pilipovi\'c, N. Teofanov} \emph{On a symbol
class of Elliptic Pseudodifferential Operators}, {Bull. Acad.
Serbe Sci. Arts}, \textbf {27} (2002), 57--68.

\bibitem{Pilipovic3} \bysame \emph{Pseudodifferential operators on
ultra-modulation spaces}, J. Funct. Anal., \textbf{208} (2004),
194--228.

\bibitem{PTT1} {S. Pilipovi\'c, N. Teofanov, J. Toft},
\emph{Micro-local analysis in Fourier Lebesgue and modulation
spaces. Part I}, preprint, in  arXiv:0804.1730, 2008.

\bibitem{PTT2} {S. Pilipovi\'c, N. Teofanov, J. Toft},
\emph{Micro-local analysis in Fourier Lebesgue and modulation
spaces. Part II}, preprint, in arXiv:0805.4476, 2008.

\bibitem{RSTT} M. Ruzhansky, m. Sugimoto, N. Tomita, J. Toft
\emph{Changes of variables in modulation and Wiener amalgam spaces},
Preprint, 2008, Available at arXiv:0803.3485v1.

\bibitem{Sjostrand1}  {J. Sj{\"o}strand} \emph{An algebra of
pseudodifferential operators}, {Math. Res. L.}, \textbf 1 (1994),
185--192.

\bibitem{Sjostrand2} \bysame \emph{Wiener type algebras of
pseudodifferential operators}, S\'eminaire Equations aux D\'eriv\'ees
Partielles, Ecole Polytechnique, 1994/1995, {Expos\'e n$^{\circ}$ IV.}

\bibitem{Sugimoto1} {M. Sugimoto, N. Tomita} \emph{The dilation
property of modulation spaces and their inclusion relation with Besov
Spaces}, {J. Funct. Anal. (1)}, \textbf{248} (2007),
79--106.

\bibitem{Tachizawa1}  {K. Tachizawa} \emph{The boundedness of
pseudo-differential operators on modulation spaces},
Math. Nachr., \textbf{168} (1994), 263--277.

\bibitem{Teofanov1}  {N. Teofanov} \emph{Ultramodulation spaces and
pseudodifferential operators}, {Endowment Andrejevi\'c}, Beograd,
2003.

\bibitem{Teofanov2} \bysame \emph{Modulation spaces, Gelfand-Shilov
spaces and pseudodifferential operators}, Sampl. Theory Signal
Image Process, \textbf{5} (2006), 225--242.

\bibitem{Toft2} J. Toft \emph{Continuity properties for
modulation spaces with applications to pseudo-differential calculus,
I}, {J. Funct. Anal.}, \textbf{207} (2004),
399--429.

\bibitem{Toft35} \bysame \emph{Convolution and embeddings for
weighted modulation spaces {\rm {in: P. Boggiatto, R. Ashino,
M. W. Wong (Eds)}} Advances in Pseudo-Differential Operators,}
Operator Theory: Advances and Applications \textbf{155},
Birkh{\"a}user Verlag, Basel 2004, pp. 165--186.

\bibitem{To8} \bysame \emph{Continuity
properties for modulation spaces with applications to
pseudo-differential calculus, II}, {Ann. Global Anal. Geom.},
\textbf{26} (2004), 73--106.

\bibitem{To9} \bysame \emph{Continuity and Schatten-von Neumann
Properties for Pseudo-Differential Operators and Toeplitz
operators on Modulation Spaces}, The Erwin Schr{\"o}dinger
International Institute for Mathematical Physics, Preprint ESI
\textbf{1732} (2005).

\bibitem{Toft4} \bysame \emph{Continuity and Schatten
properties for pseudo-differential operators on modulation spaces {\rm
{in: J. Toft, M. W. Wong, H. Zhu (Eds) Modern Trends in
Pseudo-Differential Operators,}}} Operator Theory: Advances and
Applications \textbf{172}, Birkh{\"a}user Verlag, Basel, 2007,
pp. 173--206.

\bibitem{Wong} M. W. Wong \newblock \textit{An Introduction To
Pseudodifferential Operators} 2nd Edition, World Scientific, 1999.
\end{thebibliography}
\end{document}